\documentclass[11pt]{article}
\usepackage{amscd}
\usepackage{bbm}
\usepackage{amsmath, multicol}
\usepackage{amssymb}
\usepackage{epic}
\usepackage{caption}
\usepackage{graphicx}
\usepackage{float}
\newcommand{\s}{{\lfloor s \rfloor}}
\renewcommand{\S}{{\lceil s \rceil}}
\renewcommand{\P}{{\cal P}}
\newtheorem{lm}{Lemma}[section]
\newtheorem{thm}{Theorem}[section]

\newtheorem{cor}{Corollary}[section]

\begin{document}
\title{\bf Graphs with Diameter $n-e$ \\ Minimizing the Spectral Radius}
\author{Jingfen Lan\thanks{Tsinghua University, Beijing, China, ({\tt lanjf07@mails.tsinghua.edu.cn}).}
\and
Linyuan Lu
\thanks{University of South Carolina, Columbia, SC 29208,
({\tt lu@math.sc.edu}). This author was supported in part by NSF
grant DMS 1000475. }
\and Lingsheng Shi \thanks{Tsinghua University, Beijing, China, ({\tt lshi@math.tsinghua.edu.cn}). This author was supported in part by NSFC grant 10701046.}
}
\maketitle
\begin{abstract}
  The spectral radius $\rho(G)$ of a graph $G$ is the largest
  eigenvalue of its adjacency matrix $A(G)$. For a fixed integer
  $e\ge 1$, let $G^{min}_{n,n-e}$ be a graph with minimal spectral
  radius among all connected graphs on $n$ vertices with diameter $n-e$. Let
  $P_{n_1,n_2,...,n_t,p}^{m_1,m_2,...,m_t}$ be a tree obtained from a
  path of $p$ vertices
  ($0 \sim 1 \sim 2 \sim \cdots \sim (p-1)$) by linking one pendant path $P_{n_i}$
   at $m_i$ for each $i\in\{1,2,...,t\}$. For $e=1,2,3,4,5$,
  $G^{min}_{n,n-e}$ were determined in the literature.
  Cioab\v{a}-van
  Dam-Koolen-Lee \cite{CDK} conjectured for fixed $e\geq 6$,
  $G^{min}_{n,n-e}$ is in the family
  $\P_{n,e}=\{P_{2,1,...1,2,n-e+1}^{2,m_2,...,m_{e-4},n-e-2}\mid
  2<m_2<\cdots<m_{e-4}<n-e-2\}$. For $e=6,7$, they conjectured
  $G^{min}_{n,n-6}=P^{2,\lceil\frac{D-1}{2}\rceil,D-2}_{2,1,2,n-5}$
  and $G^{min}_{n,n-7}=P^{2,\lfloor\frac{D+2}{3}\rfloor,D-
    \lfloor\frac{D+2}{3}\rfloor , D-2}_{2,1,1,2,n-6}$. In this paper,
  we settle their conjectures positively.  Note that any tree in $\P_{n,e}$
  is uniquely determined by its internal path lengths.  For any $e-4$
  non-negative integers $k_1,k_2,\ldots, k_{e-4}$, let
  $T_{(k_1,k_2,\ldots,k_{e-4})}=
  P_{2,1,...1,2,n-e+1}^{2,m_2,...,m_{e-4},n-e-2}$ with
  $k_i=m_{i+1}-m_{i}-1$, for $1\leq i \leq {e-4}$. (Here we assume $m_1=2$
and $m_{e-3}=n-e-2$.)

  Let $s= \frac{\sum_{i=1}^{e-4}k_i+2}{{e-4}}$.  For any integer $e\geq 6$ and
  sufficiently large $n$, we proved that $G^{min}_{n,n-e}$ must be one
  of the trees $T_{(k_1,k_2,\ldots,k_{e-4})}$ with the parameters
  satisfying $\s -1\leq k_j\leq \s \leq k_i \leq \S+1$ for $j=1,{e-4}$ and
  $i=2,\ldots, e-5$. Moreover, $0\le k_i-k_j\le 2$ for $2\le
  i\le{e-5}, \quad j=1,{e-4};$ and $|k_i-k_j|\le 1$ for $2\le
  i,j\le{e-5}$. These results are best possible as shown by cases $e=6,7,8$,
where  $G^{min}_{n,n-e}$ are completely determined here.
Moreover, if $n-6$ is divisible by $e-4$ and $n$ is sufficiently large, then
$G^{min}_{n,e}=T_{(k-1,k,k,\ldots,k,k,k-1)}$ where $k=\frac{n-6}{e-4}-2$.
\end{abstract}

\section{Introduction}
Let $G=(V,E)$ be a simple connected graph, and  $A(G)$ be the adjacency
matrix of $G$. The {\em characteristic polynomial} of $G$ is defined
by $\phi_G(\lambda)=\det(\lambda I-A(G))$.
The {\em spectral radius}, denoted by $\rho(G)$, is
the largest root of $\phi_G$.
The problem of determining graphs with small spectral radius can be
traced back to Hoffman and Smith \cite{smith, HS, Hoffman}. They completely
determined all connected graphs $G$ with $\rho(G)\leq 2$.  The connected
 graphs with $\rho(G)<2$ are
precisely simple Dynkin Diagrams  $A_n$, $D_n$, $E_6$, $E_7$, and $E_8$.
The connected graphs with $\rho(G)=2$ are exactly  those simple extended  Dynkin Diagrams
$\tilde A_n$, $\tilde D_n$, $\tilde E_6$, $\tilde E_7$, and $\tilde E_8$.
Cvetkovi\'c et al. \cite{CDG} gave a nearly complete description
of all graphs $G$ with $2 < \rho(G) \leq \sqrt{2 + \sqrt{5}}$. Their description
was completed by Brouwer and Neumaier \cite{BN}.
Those graphs are some special trees with at most two vertices of degree 3.
Wang et al. \cite{wang} studied some graphs with spectral radius close to $\frac{3}{2}{\sqrt{2}}$.
Woo and Neumaier \cite{WN} determined the structures of  graphs $G$ with $\sqrt{2 + \sqrt{5}}\leq \rho(G)\leq \frac{3}{2}{\sqrt{2}}$;
if  $G$ has maximum degree at least $4$, then $G$ is a {\it dagger} (i.e., a path is
attached to a leaf of a star $S_4$); if $G$ is a tree with maximum degree at most $3$, then $G$ is an {\it open quipus} (i.e., the vertices of degree $3$ lies on a path); else $G$ is a {\it closed quipus} (i.e., a unicyclic graph with maximum degree at most $3$ satisfies that  the vertices of degree $3$ lies on a cycle).

\begin{figure}[htpb]
\begin{center}
\setlength{\unitlength}{0.7cm}
\begin{picture}(12,4)
\multiput(0,1)(4,0){4}{\circle*{0.2}}
\multiput(1,1)(2,0){6}{\circle*{0.2}}
\multiput(4,2)(4,0){2}{\circle*{0.2}}
\multiput(4,4)(4,0){2}{\circle*{0.2}}
\multiput(0,1)(11,0){2}{\line(1,0){1}}
\multiput(4,1)(4,0){2}{\line(0,1){1}}
\multiput(3,1)(4,0){2}{\line(1,0){2}}
\dashline{0.3}(4,2)(4,4)\dashline{0.3}(8,2)(8,4)\dashline{0.3}(1,1)(11,1)
\multiput(5,3)(1,0){3}{\circle*{0.1}}
\put(0,0.2){$0$}\put(1,0.2){$1$}\put(11.7,0.2){$p-1$}
\put(4,0.2){$m_1$}\put(8,0.2){$m_t$}
\put(3,3){$P_{n_1}$}\put(8.3,3){\bf $P_{n_t}$}
\end{picture}
\caption{$P_{n_1,n_2,...,n_t,p}^{m_1,m_2,...,m_t}$}\label{p1}
\end{center}
\end{figure}
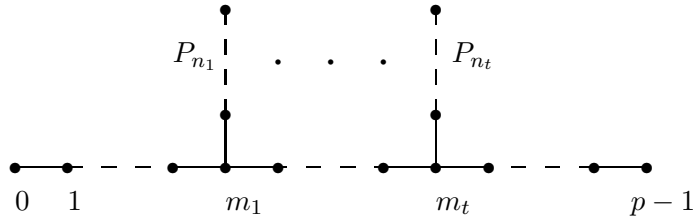

Van Dam and Kooij \cite{DK} used the following notation to denote an
open quipus.  Let $P_{n_1,n_2,...n_t,p}^{m_1,m_2,...,m_t}$ be a tree obtained from
 a path on $p$ vertices ($0 \sim 1 \sim 2 \sim \cdots \sim (p-1)$) by linking one pendant path $P_{n_i}$
at $m_i$ for $i=1,2,...,t$ (see
Figure \ref{p1}.)  The path $0 \sim 1 \sim 2 \sim \cdots \sim (p-1)$ is called {\em main path}.
For $i=1,...,t-1$, let $P^{(i)}$ be the $i$-th
{\it internal path} ($m_i\sim m_i\!+\!1 \sim \cdots \sim m_{i+1}$) and $k_i=m_{i+1}-m_i-1$ be the
the number of internal vertices on $P^{(i)}$. In general, an {\em
  internal path} in $G$ is a path
$v_0\sim v_1\sim \cdots\sim v_s$ such that $d(v_0)>2$, $d(v_s)>2$, and
$d(v_i)=2$, whenever $0<i<s$.  An internal path is {\em closed} if
$v_0=v_s$.

Recently van Dam and Kooij \cite{DK} asked an interesting question {\em ``which connected graph of order $n$ with
a given diameter $D$ has minimal spectral radius}?''.
Here  the {\em diameter} of a connected graph is the
maximum distance among all pairs of its vertices.
 They \cite{DK} solved this
problem explicitly for graphs with diameter $D\in\{1, 2, \lfloor
n/2\rfloor, n-3, n-2, n-1\}$.   The cases $D=1$ and $D=n-1$ are trivial.
A {\em minimizer} graph, denoted by $G_{n,D}^{min}$, is a graph
that has the minimal spectral radius among
all the graphs of order $n$ and diameter $D$.
They proved that $G_{n,2}^{min}$ is either a star or a Moore graph;
$G_{n,\lfloor n/2\rfloor}^{min}$  is the cycle $C_n$;
$G_{n,n-2}^{min}$ is the tree $P_{1,n-1}^{1}$;
$G_{n,n-3}^{min}$ is the tree $P_{1,1,n-2}^{1,n-4}$.
They conjectured
$G^{min}_{n,n-e}=P^{\lfloor \frac{e-1}{2}\rfloor, n-e-\lceil\frac{e-1}{2}\rceil}
_{\lfloor \frac{e-1}{2}\rfloor, \lceil\frac{e-1}{2}\rceil,n-e+1}$  for any constant $e\geq 1$
and $n$ large enough.

This conjecture is proved for  $e=4$ by Yuan et al. \cite{YSL} and for  $e=5$ by Cioab\v{a} et al.
\cite{CDK}. However, it is disproved for $e=6$ by \cite{Sun}
and  any $e\geq 6$ by \cite{CDK} when $n$ large enough.
Cioab\v{a}-van Dam-Koolen-Lee \cite{CDK} proved the following theorem.

{\bf Theorem 5.2 of \cite{CDK}:} \label{page2}
{\it For $e\ge6$,
$\rho(G_{n,n-e}^{min}) \to \sqrt{2+\sqrt 5}$ as $n \to \infty$.
Moreover
$G_{n,n-e}^{min}$ is contained in one of the following three families of graphs
\begin{eqnarray*}
&&\P_{n,e}=\left\{P_{2,1,...1,2,n-e+1}^{2,m_2,...,m_{e-4},n-e-2}\mid 2<m_2<\cdots<m_{e-4}<n-e-2\right\},\\
&&\P'_{n,e}=\left\{P_{2,1,...1,1,n-e+1}^{2,m_2,...,m_{e-3},n-e-1}\mid 2<m_2<\cdots<m_{e-3}<n-e-1\right\},\\
&&\P''_{n,e}=\left\{P_{1,1,...1,1,n-e+1}^{1,m_2,...,m_{e-2},n-e-1}\mid 1<m_2<\cdots<m_{e-2}<n-e-1\right\}.
\end{eqnarray*}
}
Cioab\v{a} et al. \cite{CDK} made three conjectures.
\begin{description}
\item[Conjecture 1] (\cite{CDK}, 5.3) For fixed $e\geq 5$, a minimizer
  graph with $n$ vertices and diameter $D=n-e$ is in the family $\P_{n,e}$,
  for $n$ large enough.
\item [Conjecture 2] (\cite{CDK}, 5.4)  The graph $P^{2,\lceil\frac{D-1}{2}\rceil,D-2}_{2,1,2,n-5}$
 is the unique minimizer graph with $n$ vertices and diameter
$D = n-6$, for $n$ large enough.
\item [Conjecture 3] (\cite{CDK}, 5.5)  The graph
$P^{2,\lfloor\frac{D+2}{3}\rfloor,D- \lfloor\frac{D+2}{3}\rfloor
, D-2}_{2,1,1,2,n-6}$
 is the unique minimizer graph with $n$ vertices and diameter
$D = n-7$, for $n$ large enough.
\footnote{Conjecture 5.5 of \cite{CDK} contains a typo:
``...$P^{2,\lfloor\frac{D-2}{3}\rfloor,D- \lfloor\frac{D-2}{3}\rfloor
, D-2}_{2,1,1,2,n-6}$...''.}
\end{description}

In this paper, we settle these three conjectures positively.

\begin{figure}[hbt]
\begin{center}
\setlength{\unitlength}{.8cm}
\begin{picture}(17,9)
\multiput(0,7)(1,0){3}{\circle*{0.18}}
\multiput(11,7)(1,0){3}{\circle*{0.18}}
\multiput(2,8)(0,1){2}{\circle*{0.18}}
\multiput(11,8)(0,1){2}{\circle*{0.18}}
\multiput(4,7)(0,1){2}{\circle*{0.18}}
\multiput(6,7)(0,1){2}{\circle*{0.18}}
\multiput(9,7)(0,1){2}{\circle*{0.18}}
\multiput(0,7)(11,0){2}{\line(1,0){2}}
\multiput(2,7)(9,0){2}{\line(0,1){2}}
\multiput(4,7)(2,0){2}{\line(0,1){1}}
\multiput(9,7)(0,0){1}{\line(0,1){1}}
\dashline{0.2}(2,7)(11,7)
\multiput(7,7.5)(0.3,0){3}{\circle*{0.1}}
\multiput(7,6.4)(0.3,0){3}{\circle*{0.1}}
\put(2.8,6.2){$k_1$}\put(4.85,6.2){$k_2$}\put(9.9,6.2){$k_{e-4}$}
\put(2.3,6.9){$\underbrace{\quad\quad\quad}$}
\put(4.32,6.9){$\underbrace{\quad\quad\quad}$}
\put(9.35,6.9){$\underbrace{\quad\quad\quad}$}
\put(15,7){$T_{(k_1,k_2,...,k_{e-4})}$}
\multiput(0,3.5)(1,0){3}{\circle*{0.18}}
\multiput(11,3.5)(1,0){2}{\circle*{0.18}}
\multiput(2,4.5)(0,1){2}{\circle*{0.18}}
\multiput(11,4.5)(0,0){1}{\circle*{0.18}}
\multiput(4,3.5)(0,1){2}{\circle*{0.18}}
\multiput(6,3.5)(0,1){2}{\circle*{0.18}}
\multiput(9,3.5)(0,1){2}{\circle*{0.18}}
\multiput(0,3.5)(0,0){1}{\line(1,0){2}}
\multiput(11,3.5)(0,0){1}{\line(1,0){1}}
\multiput(2,3.5)(0,0){1}{\line(0,1){2}}
\multiput(4,3.5)(2,0){2}{\line(0,1){1}}
\multiput(9,3.5)(2,0){2}{\line(0,1){1}}
\dashline{0.2}(2,3.5)(11,3.5)
\multiput(7,4)(0.3,0){3}{\circle*{0.1}}
\multiput(7,2.9)(0.3,0){3}{\circle*{0.1}}
\put(2.8,2.7){$k_1$}\put(4.85,2.7){$k_2$}\put(9.9,2.7){$k_{e-3}$}
\put(2.3,3.4){$\underbrace{\quad\quad\quad}$}
\put(4.32,3.4){$\underbrace{\quad\quad\quad}$}
\put(9.35,3.4){$\underbrace{\quad\quad\quad}$}
\put(15,3.5){$T'_{(k_1,k_2,...,k_{e-3})}$}
\multiput(1,1)(1,0){2}{\circle*{0.18}}
\multiput(11,1)(1,0){2}{\circle*{0.18}}
\multiput(2,2)(9,0){2}{\circle*{0.18}}
\multiput(4,1)(0,1){2}{\circle*{0.18}}
\multiput(6,1)(0,1){2}{\circle*{0.18}}
\multiput(9,1)(0,1){2}{\circle*{0.18}}
\multiput(1,1)(10,0){2}{\line(1,0){1}}
\multiput(2,1)(9,0){2}{\line(0,1){1}}
\multiput(4,1)(2,0){2}{\line(0,1){1}}
\multiput(9,1)(0,0){1}{\line(0,1){1}}
\dashline{0.2}(2,1)(11,1)
\multiput(7,1.5)(0.3,0){3}{\circle*{0.1}}
\multiput(7,0.4)(0.3,0){3}{\circle*{0.1}}
\put(2.8,.2){$k_1$}\put(4.85,.2){$k_2$}\put(9.9,.2){$k_{e-2}$}
\put(2.3,0.9){$\underbrace{\quad\quad\quad}$}
\put(4.32,0.9){$\underbrace{\quad\quad\quad}$}
\put(9.35,0.9){$\underbrace{\quad\quad\quad}$}
\put(15,1){$T''_{(k_1,k_2,...,k_{e-2})}$}
\end{picture}
\end{center}
\caption{The three families of graphs: $\P_{n,e}$, $\P'_{n,e}$,$\P''_{n,e}$.}\label{p2}
\end{figure}
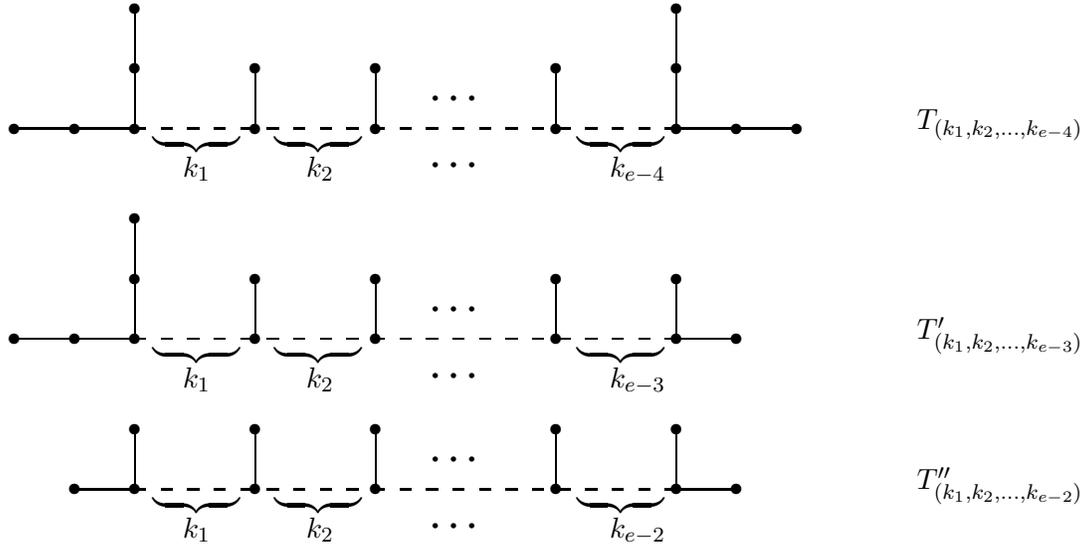

Note that graphs in each family can be determined by the lengths of
internal paths (see Figure~ \ref{p2}). The parameters
$k_i$'s and $m_i$'s are related as follows.
In the first family $\P_{n,e}$,
$T_{(k_1,k_2,...,k_{e-4})}= P_{2,1,...1,2,n-e+1}^{2,m_2,...,m_{e-4},n-e-2}$
if $k_i=m_{i+1}-m_i-1$ for $1\leq i\leq e-4$, where
 $m_1=2$ and $m_{e-3}=n-e-2$. In the second family $\P'_{n,e}$,
$T'_{(k_1,k_2,...,k_{e-3})}=P_{2,1,...1,1,n-e+1}^{2,m_2,...,m_{e-3},n-e-1}$
if $k_i=m_{i+1}-m_i-1$ for $1\leq i\leq e-3$,
where  $m_1=2$ and $m_{e-2}=n-e-1$. In the third family $\P''_{n,e}$,
$T''_{(k_1,k_2,...,k_{e-2})}=P_{1,1,...1,1,n-e+1}^{1,m_2,...,m_{e-2},n-e-1}$
if $k_i=m_{i+1}-m_i-1$ for $1\leq i\leq e-2$,
where $m_1=1$ and $m_{e-1}=n-e-1$. In all three cases, 
the summation of all $k_i$'s is always equal to $n-2e$.


\medskip
We have the following theorem.

\begin{thm}\label{thmMain}
For any $e\ge 6$ and sufficiently large $n$,
 $G^{min}_{n, n-e}$ must be  a tree $T_{ (k_1,k_2,\ldots,k_{e-4})}$ in $\P_{n,e}$ satisfying
\begin{enumerate}
\item $\s -1 \leq k_j\leq \s \leq k_i\leq \S+1$ for $2\leq i\leq e-5$
and $j=1, e-4$, where $s= \frac{n-6}{e-4}-2$,
\item $0\le k_i-k_j\le 2$ for $2\le i\le{e-5}$ and $j=1,e-4$,
\item $|k_i-k_j|\le 1$ for $2\le i, j\le{e-5}$.
\end{enumerate}
In particular, if $n-6$ is divisible by $e-4$, then $G^{min}_{n, n-e}=T_{ (s-1,s,
\ldots, s, s-1)}$.
\end{thm}

Here we  completely determine the $G_{n,n-e}^{min}$ for $e=6,7,8$ and settle the conjectures
 2 and 3 positively.

\begin{thm}\label{thm1.2}
For  $e=6$ and $n$ large enough, $G_{n,n-e}^{min}$ is
unique up to  a graph isomorphism.
\begin{enumerate}
\item If $n=2k+12$, then $G_{n,n-6}^{min}=T_{(k,k)}$.
\item  If $n=2k+13$, then $G_{n,n-6}^{min}=T_{(k,k+1)}$.
\end{enumerate}
\end{thm}

\begin{thm}\label{thm1.3}
  For  $e=7$ and $n$ large enough, $G_{n,n-e}^{min}$ is
unique up to  a graph isomorphism.
\begin{enumerate}
\item If $n=3k+14$, then $G_{n,n-7}^{min}=T_{(k,k,k)}$.
\item  If $n=3k+15$, then $G_{n,n-7}^{min}=T_{(k,k+1,k)}$.
 \item If $n=3k+16$, then $G_{n,n-7}^{min}=T_{(k,k+2,k)}$.
\end{enumerate}
\end{thm}

\begin{thm}\label{thm1.4}
  For $e=8$ and $n$ large enough, $G_{n,n-e}^{min}$ is determined
up to  a graph isomorphism as follows.
\begin{enumerate}
\item If $n=3k+16$, then $G_{n,n-8}^{min}=T_{(k,k,k,k)}$, $T_{(k,k,k+1,k-1)}$, or $T_{(k-1,k+1,k+1,k-1)}$;
all three trees have the same spectral radius.
\item  If $n=3k+17$, then $G_{n,n-8}^{min}=T_{(k,k+1,k,k)}$.
 \item If $n=3k+18$, then $G_{n,n-8}^{min}=T_{(k,k+1,k+1,k)}$.
 \item If $n=3k+19$, then $G_{n,n-8}^{min}=T_{(k,k+1,k+2,k)}$.
\end{enumerate}
\end{thm}

For $e=6$, Theorem \ref{thm1.2} is an easy corollary of
Theorem \ref{thmMain}.   Theorem \ref{thm1.3} and Theorem
\ref{thm1.4} show that the bounds on $k_i$'s in Theorem \ref{thmMain}
are best possible.

The remaining of the paper is organized as follows. In section 2, we  prove
some useful lemmas. The proof of Theorem \ref{thmMain} is presented
in section 3 and the proofs of Theorem \ref{thm1.3} and \ref{thm1.4} are given
in section 4.

\section{Basic notations and Lemmas}
\subsection{Preliminary results}
For any vertex $v$ in a graph $G$, let $N(v)$ be the neighborhood of $v$.
Let $G-v$ be the remaining graph of $G$ after deleting the vertex $v$
(and all edges incident to $v$).
Similarly, $G-u-v$ is the remaining graph of $G$ after deleting two vertices $u,v$.
Here are some basic facts found in literature \cite{CDS, HS, parlett, Sun}, which will be used later.
\begin{lm}{\em\cite{CDS}}\label{lm2.1}
Suppose that $G$ is a connected graph. If $v$ is not in any cycle of $G$, then
$\phi_G=\lambda\phi_{G-v}-\sum_{w\in N(v)}\phi_{G-w-v}.$
If $e=uv$ is a cut edge of $G$, then
$\phi_G=\phi_{G-e}-\phi_{G-u-v}.$
\end{lm}

\begin{lm}{\em\cite{CDS}}\label{lm2.2} Let $G_1$ and $G_2$ be two
graphs, then the following statements hold.
\begin{enumerate}
\item If $G_2$ is a proper subgraph of $G_1$, then $\rho(G_1)>\rho(G_2)$.
\item If $G_2$ is a spanning proper subgraph of $G_1$, then  $\rho(G_1)>\rho(G_2)$ and $\phi_{G_2}(\lambda)>\phi_{G_1}(\lambda)$ for all $\lambda\ge\rho(G_1)$.
\item If $\phi_{G_2}(\lambda)>\phi_{G_1}(\lambda)$ for all $\lambda\ge{\rho(G_1)}$, then $\rho(G_2)<\rho(G_1)$.
\item If $\phi_{G_1}(\rho(G_2))<0$,then $\rho(G_1)>\rho(G_2)$.
\end{enumerate}
\end{lm}

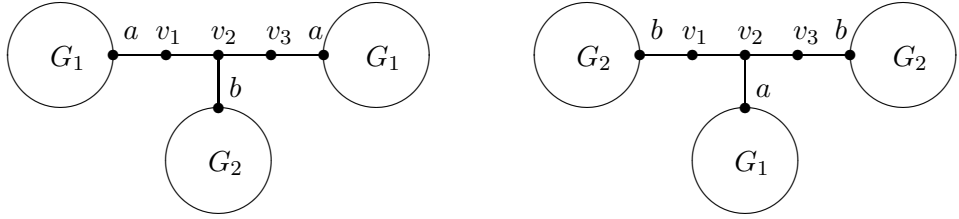
\begin{figure}[h]
\begin{center}
\setlength{\unitlength}{0.7cm}
\begin{picture}(20,3)
\multiput(5,0)(3,0){1}{\circle{2}}
\multiput(2,2)(6,0){2}{\circle{2}}
\multiput(5,1)(0,1){1}{\circle*{0.2}}
\multiput(3,2)(1,0){5}{\circle*{0.2}}
\multiput(5,1)(1,0){1}{\line(0,1){1}}
\multiput(3,2)(1,0){1}{\line(1,0){4}}
\put(4.8,-0.2){$G_2$}\put(5.2,1.2){$b$}\put(3.8,2.3){$v_1\quad
v_2\quad v_3$}
\multiput(3.2,2.3)(3.5,0){2}{$a$}
\multiput(1.8,1.8)(6,0){2}{$G_1$}
\multiput(15,0)(3,0){1}{\circle{2}}
\multiput(12,2)(6,0){2}{\circle{2}}
\multiput(15,1)(0,1){1}{\circle*{0.2}}
\multiput(13,2)(1,0){5}{\circle*{0.2}}
\multiput(15,1)(1,0){1}{\line(0,1){1}}
\multiput(13,2)(1,0){1}{\line(1,0){4}}
\put(14.8,-0.2){$G_1$}\put(15.2,1.2){$a$}\put(13.8,2.3){$v_1\quad
v_2\quad v_3$}
\multiput(13.2,2.3)(3.5,0){2}{$b$}
\multiput(11.8,1.8)(6,0){2}{$G_2$}
\end{picture}
\end{center}
\caption{The graphs $H_1$ and $H_2$}\label{p3}
\end{figure}

\begin{lm}\label{lm2.3}{\rm\cite{Sun}} Let $G_1$ and $G_2$ be two (possibly empty) graphs with $a\in V(G_1)$ and $b\in V(G_2)$, and let $H_1$ and $H_2$ be two graphs shown in Figure \ref{p3}. Then $\rho(H_1)=\rho(H_2)$.
\end{lm}

\begin{lm}{\em\cite{HS}}\label{lm2.4} Let $uv$ be an edge of a connected graph $G$ of order $n$, and denote by $G_{u,\,v}$
the graph obtained from $G$ by subdividing the edge $uv$ once, i.e.,
adding a new vertex $w$ and edges $wu,wv$ in $G-uv$. Then the
following two properties hold.
\begin{enumerate}
\item If $uv$ does not belong to an internal path of $G$ and $G\neq
C_n$, then $\rho(G_{u,\,v})>\rho(G)$.
\item If $uv$ belongs to an internal path of $G$ and $G\neq P_{1,1,n}^{1,n-2}$, then $\rho(G_{u,\,v})<\rho(G)$.
\end{enumerate}
\end{lm}

\begin{thm}[Cauchy Interlace Theorem (see p.183, \cite{parlett})]
Let $A$ be a Hermitian matrix of order $n$,
and let $B$ be a principal submatrix of $A$ of order $n-1$.
If $\lambda_n\leq \lambda_{n-1}\leq \cdots\leq \lambda_1$
lists the eigenvalues of $A$ and  $\mu_{n-1}\leq \mu_{n-2}\leq \cdots\leq \mu_1$
the eigenvalues of $B$, then
$$\lambda_n\leq \mu_{n-1}\leq \lambda_{n-1}\leq \cdots\leq \lambda_2 \leq \mu_1 \leq \lambda_1.$$
\end{thm}

Applying Cauchy Interlace Theorem to the adjacency matrices of graphs,
we have the following corollary.
\begin{cor}
  \label{cor0}
Suppose $G$ is a connected graph. Let $\lambda_2(G)$ be the second
largest eigenvalue of $G$. For any vertex $v$, we have
$$\lambda_2(G)<\rho(G-v)<\rho(G).$$
\end{cor}

\subsection{Our approach}
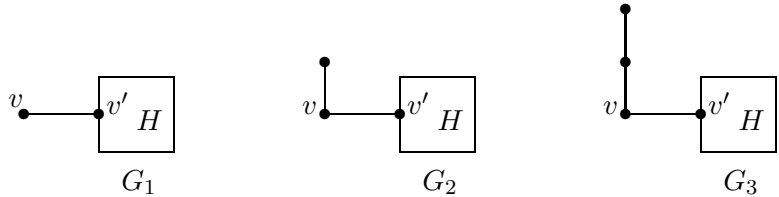
\begin{figure}[htbp]
\begin{center}
\setlength{\unitlength}{1.0cm}
\begin{picture}(12,2.5)
\multiput(1,1)(4,0){3}{\circle*{0.15}}
\multiput(2,1)(4,0){3}{\circle*{0.15}}
\multiput(5,1.7)(4,0){2}{\circle*{0.15}}
\multiput(9,2.4)(0,0){1}{\circle*{0.15}}
\multiput(1,1)(4,0){3}{\line(1,0){1}}
\multiput(2,0.5)(4,0){3}{\line(1,0){1}}
\multiput(2,1.5)(4,0){3}{\line(1,0){1}}
\multiput(2,0.5)(4,0){3}{\line(0,1){1}}
\multiput(3,0.5)(4,0){3}{\line(0,1){1}}
\multiput(5,1)(0,0){1}{\line(0,1){0.7}}
\multiput(9,1)(0,0){1}{\line(0,1){1.4}}
\put(0.8,1.1){$v$}\put(4.7,1){$v$}\put(8.7,1){$v$}
\put(2.1,1){$v'$}\put(6.1,1){$v'$}\put(10.1,1){$v'$}
\put(2.5,0.8){$H$}\put(6.5,0.8){$H$}\put(10.5,0.8){$H$}
\put(2.3,0){$G_1$}\put(6.3,0){$G_2$}\put(10.3,0){$G_3$}
\end{picture}
\end{center}
\caption{For $i=1,2,3$, three graphs $(G_i, v)$ are constructed from $(H,v')$.}\label{p4}
\end{figure}

A {\em rooted graph} $(G,v)$ is a graph $G$ together with a designated vertex $v$ as a root.
For $i=1,2,3$ and a given rooted graph $(H,v')$,
we  get a new rooted graph $(G_i,v)$  from $H$ by attaching a path $P_i$ to  $v'$ and
changing the root from $v'$ to $v$ as shown by Figure \ref{p4}.

Note that any tree in the three families $\P_{n,e}$,
$\P'_{n,e}$, $\P''_{n,e}$ can be built up from a single vertex through
a sequence of three operations above.  Applying Lemma \ref{lm2.1}, we
observe that the pair $(\phi_{G_i}, \phi_{G_i-v})$ linearly depends on
$(\phi_H, \phi_{H-v'})$ with coefficients in
$\mathbb{Z}[\lambda]$. We can choose proper base to diagonalize
the operation from $(H,v')$ to
$(G_1,v)$.

Let $\lambda_0$ be the constant $\sqrt{2+\sqrt 5}=2.058\cdots$.
In this paper, we consider only the range
$\lambda\geq \lambda_0$.
Let $x_1$ and $x_2$ be two roots of the equation $x^2-\lambda x +1=0$. We have
 $$x_1=\frac{\lambda-\sqrt{\lambda^2-4}}{2}, \quad x_2=\frac{\lambda+\sqrt{\lambda^2-4}}{2}$$
and
\begin{equation}\label{eq1}
x_1+x_2=\lambda, \quad x_1x_2=1.\\
\end{equation}
For any vertex $v$ in a graph $G$, we define two functions (of $\lambda$) $p_{(G,v)}$ and $q_{(G,v)}$ satisfying
\begin{eqnarray*}
   \phi_G&=&p_{(G,v)}+q_{(G,v)}, \\
  \phi_{G-v}&=&x_2p_{(G,v)} + x_1q_{(G,v)}.
\end{eqnarray*}
This definition can be written in the following matrix form:
\begin{equation}\label{eq2}
\left(
  \begin{array}[c]{c}
    \phi_G\\
   \phi_{G-v}
  \end{array}
\right)
=\left (\begin{array}{ll}
1 & 1\\ x_2 & x_1
\end{array}\right )
\left(
  \begin{array}[c]{c}
    p_{(G,v)}\\
   q_{(G,v)}
  \end{array}
\right).
\end{equation}

Using Equation (\ref{eq1}),
 we can solve $p_{(G,v)}$ and $q_{(G,v)}$ and get
\begin{equation}\label{eq3}
\left(\begin{array}[c]{c} p_{(G,v)}\\ q_{(G,v)}\end{array}\right)=\frac{1}{x_2-x_1}\left (\begin{array}{ll} -x_1 & 1\\ x_2 & -1 \end{array}\right )
\left(
  \begin{array}[c]{c}
    \phi_G\\
   \phi_{G-v}
  \end{array}
\right).
\end{equation}

For example, let $v$ be the center of the odd path $P_{2k+1}$.
We have

\begin{eqnarray}
\label{eq:P1}
\left(\!\!
  \begin{array}[c]{c}
    p_{(P_1,v)}\\
   q_{(P_1,v)}
  \end{array}\!\!
\right)
&=&
\frac{1}{x_2-x_1}
\left(\!\!
\begin{array}[c]{c}
   -x_1^2\\
  x_2^2
  \end{array}\!\!
\right),\\
\label{eq:P3}
\left(\!\!
  \begin{array}[c]{c}
    p_{(P_3,v)}\\
   q_{(P_3,v)}
  \end{array}\!\!
\right)
&=&
\lambda
\left(\!\!
\begin{array}[c]{c}
   x_1^2\\
  x_2^2
  \end{array}\!\!
\right),\\
\label{eq:P5}
\left(\!\!
  \begin{array}[c]{c}
    p_{(P_5,v)}\\
   q_{(P_5,v)}
  \end{array}\!\!
\right)
&=&
\frac{\lambda^2-1}{x_2-x_1}
\left(\!\!
\begin{array}[c]{c}
   (\lambda-x_1^3 )x_1\\
 (x_2^3-\lambda) x_2
  \end{array}\!\!
\right).
\end{eqnarray}

We have the following lemma.
\begin{lm} \label{lmq}
For any tree $G$ and any vertex $v$, we have
\begin{equation}
  \label{eq:infty}
  \lim_{\lambda\to +\infty}q_{(G,v)}(\lambda)=+\infty.
\end{equation}
\end{lm}
{\bf Proof }
From Lemma \ref{lm2.1}, we have
$$\phi_G=\lambda\phi_{G-v}-\sum_{w\in N(v)}\phi_{G-w-v}.$$
By Equation (\ref{eq3}), we get
\begin{eqnarray*}
q_{(G,v)}&=&\frac{1}{x_2-x_1} (x_2  \phi_G -  \phi_{G-v})\\
 &=& \frac{1}{x_2-x_1}\left (x_2\left (\lambda\phi_{G-v}-\sum_{w\in N(v)}\phi_{G-w-v}\right ) - \phi_{G-v}\right)\\
&=& \frac{1}{x_2-x_1}\left ((\lambda x_2-1)\phi_{G-v} -x_2\sum_{w\in N(v)}\phi_{G-w-v}\right )\\
&=& \frac{x_2}{x_2-x_1}\left (x_2\phi_{G-v}-\sum_{w\in N(v)}\phi_{G-w-v}\right ).
\end{eqnarray*}
Note that $\phi_{G-v}$ is a polynomial of degree $n-1$ with highest coefficient $1$ while
$\phi_{G-w-v}$ is a polynomial of degree $n-2$  with highest coefficient $1$.
Since $x_2> 1 > x_1$, we have $x_2\phi_{G-v}-\sum_{w\in N(v)}\phi_{G-w-v}$ goes to infinity as $\lambda$ approaches infinity.
\hfill $\square$

\begin{lm}\label{lm2.6}
Let $G_1,G_2,G_3$ be the graphs shown in Figure \ref{p4}. Then the following equations hold.
\begin{enumerate}
\item $\left(\begin{array}[c]{c} p_{(G_1,v)}\\ q_{(G_1,v)}\end{array}\right)= \left (\begin{array}{ll} x_1 & 0\\ 0 & x_2 \end{array}\right ) \left(\begin{array}[c]{c} p_{(H,v')}\\ q_{(H,v')}\end{array}\right).\\$
\item $\left(\begin{array}[c]{c} p_{(G_2,v)}\\ q_{(G_2,v)}\end{array}\right)= \dfrac{1}{x_2-x_1}\left (\begin{array}{ll} \lambda-x_1^3 & x_1\\ -x_2 & x_2^3-\lambda \end{array}\right ) \left(\begin{array}[c]{c} p_{(H,v')}\\ q_{(H,v')}\end{array}\right).\\$
\item $\left(\begin{array}[c]{c} p_{(G_3,v)}\\ q_{(G_3,v)}\end{array}\right)= \dfrac{1}{x_2-x_1}\left (\begin{array}{ll} -x_1^4+\lambda^2-1 & \lambda x_1\\ -\lambda x_2 & x_2^4-\lambda^2+1 \end{array}\right ) \left(\begin{array}[c]{c} p_{(H,v')}\\ q_{(H,v')}\end{array}\right).$
\end{enumerate}
\end{lm}
{\bf Proof } By Lemma  \ref{lm2.1}, we have
$$\left(\begin{array}[c]{c}\phi_{G_1}\\ \phi_{G_1-v} \end{array}\right)=\left (\begin{array}{ll} \lambda & -1\\ 1 & 0 \end{array}\right )\left(\begin{array}[c]{c}\phi_H\\ \phi_{H-v'} \end{array}\right).$$
Combining it with equations ( \ref{eq2}) ( \ref{eq3}), we get
\begin{eqnarray*}
\left(\begin{array}[c]{c} p_{(G_1,v)}\\ q_{(G_1,v)}\end{array}\right)&=&\left (\begin{array}{ll} 1 & 1\\ x_2 & x_1 \end{array}\right )^{-1}\left (\begin{array}{ll} \lambda & -1\\ 1 & 0 \end{array}\right )\left (\begin{array}{ll} 1 & 1\\ x_2 & x_1 \end{array}\right ) \left(\begin{array}[c]{c} p_{(H,v')}\\ q_{(H,v')}\end{array}\right)\\
&=& \frac{1}{x_2-x_1}\left (\begin{array}{ll} 2-\lambda x_1 & x_1^2-\lambda x_1+1\\ -x_2^2+\lambda x_2-1 & \lambda x_2-2 \end{array}\right) \left(\begin{array}[c]{c} p_{(H,v')}\\ q_{(H,v')}\end{array}\right)\\
&=& \left (\begin{array}{ll} x_1 & 0\\ 0 & x_2 \end{array}\right ) \left(\begin{array}[c]{c} p_{(H,v')}\\ q_{(H,v')}\end{array}\right).
\end{eqnarray*}
The proofs of items 2 and 3 are similar as that of item 1. $\hfill\Box$

We denote the three matrices by $A$, $B$, and $C$. Namely,
$$A=\left (\begin{array}{ll} x_1 & 0\\ 0 & x_2 \end{array}\right ),
B=\frac{1}{x_2-x_1} \left(\begin{array}{ll} \lambda-x_1^3 & x_1\\ -x_2 & x_2^3-\lambda \end{array}\right),
C=\frac{1}{x_2-x_1} \left(\begin{array}{ll} -x_1^4+\lambda^2-1 & \lambda x_1\\ -\lambda x_2 & x_2^4-\lambda^2+1 \end{array}\right).$$

The diagonal elements of $B$ are very useful parameters. To simplify our notations later, we define two parameters $d_1$ and $d_2$ as follows:
\begin{eqnarray}
  \label{eq:d1}
  d_1&=&\lambda-x_1^3,\\
 \label{eq:d1}
  d_2&=&x_2^3-\lambda.
\end{eqnarray}
Note that $d_2=0$ if $\lambda=\lambda_0$. The
equation (\ref{eq:P5}) can be written as
\begin{equation}
  \label{eq:P5d}
\left(\!\!
  \begin{array}[c]{c}
    p_{(P_5,v)}\\
   q_{(P_5, v)}
  \end{array}\!\!
\right)
=
\frac{\lambda^2-1}{x_2-x_1}
\left(\!\!
\begin{array}[c]{c}
   d_1x_1\\
 d_2 x_2
  \end{array}\!\!
\right).
\end{equation}
From the definitions of $d_1$ and $d_2$,
we can derive the following identity
\begin{equation}
  \label{eq:d1d2}
d_1x_2-d_2x_1=2.
\end{equation}

Given two rooted graphs $(H_1,v_1)$ and $(H_2,v_2)$, we define some new graphs. Denote by
$(H_1,v_1)\cdot P_i$, the graph consisting of the graph $H_1$ and a path $P_i$ linking one of its ends at the vertex $v_1$.
Similarly denote by $(H_1,v_1)\cdot P_i \cdot (H_2,v_2)$ the graph consisting of graphs $H_1, H_2$ and a path $P_i$ linking the two ends at $v_1, v_2$ respectively.

\begin{figure}[h]
\begin{center}
\setlength{\unitlength}{1.0cm}
\begin{picture}(7,1.3)
\multiput(1,0.5)(1,0){3}{\circle*{0.15}}
\multiput(1,0.5)(0,0){1}{\line(1,0){2}}
\multiput(0,0)(3,0){2}{\line(1,0){1}}
\multiput(0,0)(3,0){2}{\line(0,1){1}}
\multiput(0,1)(3,0){2}{\line(1,0){1}}
\multiput(1,0)(3,0){2}{\line(0,1){1}}
\put(0.1,0.4){$H_1$}\put(3.5,0.4){$H_2$}
\put(0.65,0.3){$v_1$}\put(3.1,0.3){$v_2$}\put(2,0.6){$u$}
\end{picture}
\end{center}
\caption{The graph $(H_1,v_1)\cdot P_1 \cdot (H_2,v_2)$}\label{p5}
\end{figure}
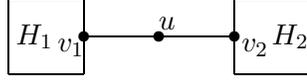

\begin{lm}\label{lm2.7} $\phi_{(H_1,v_1)\cdot P_1 \cdot (H_2,v_2)}(\lambda)=(x_2-x_1)(q_{(H_1,v_1)}q_{(H_2,v_2)}-p_{(H_1,v_1)}p_{(H_2,v_2)})$.
\end{lm}
\noindent{\bf Proof } By Lemmas \ref{lm2.1} and \ref{lm2.6}, we have
\begin{eqnarray*}
&& \hspace*{-1cm}
\phi_{(H_1,v_1)\cdot P_1 \cdot (H_2,v_2)}(\lambda) \\
&=&\lambda\phi_{H_1}\phi_{H_2}-\phi_{H_1-v_1}\phi_{H_2}-\phi_{H_2-v_2}\phi_{H_1}\\
&=&(x_1+x_2)(p_{(H_1,v_1)}+q_{(H_1,v_1)})(p_{(H_2,v_2)}+q_{(H_2,v_2)})-(p_{(H_1,v_1)}x_2+q_{(H_1,v_1)}x_1)\\
& &(p_{(H_2,v_2)}+q_{(H_2,v_2)})-(p_{(H_1,v_1)}+q_{(H_1,v_1)})(p_{(H_2,v_2)}x_2+q_{(H_2,v_2)}x_1)\\
&=&(p_{(H_2,v_2)}+q_{(H_2,v_2)})(p_{(H_1,v_1)}x_1+q_{(H_1,v_1)}x_2)-(p_{(H_1,v_1)}+q_{(H_1,v_1)})(p_{(H_2,v_2)}x_2+q_{(H_2,v_2)}x_1)\\
&=& p_{(H_2,v_2)}p_{(H_1,v_1)}x_1+q_{(H_2,v_2)}q_{(H_1,v_1)}x_2-p_{(H_2,v_2)}p_{(H_1,v_1)}x_2-q_{(H_2,v_2)}q_{(H_1,v_1)}x_1\\
&=& (x_2-x_1)(q_{(H_1,v_1)}q_{(H_2,v_2)}-p_{(H_1,v_1)}p_{(H_2,v_2)}). \hspace*{7.5cm}\square
\end{eqnarray*}

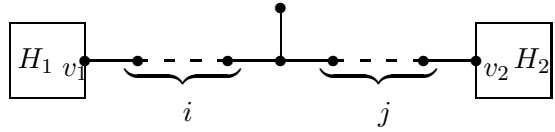
\begin{figure}[h]
\begin{center}
\setlength{\unitlength}{1.0cm}
\begin{picture}(7.5,1.6)
\multiput(1,0.5)(4.5,0){2}{\circle*{0.15}}
\multiput(1.7,0.5)(4.5,0){2}{\circle*{0.15}}
\multiput(2.9,0.5)(0.7,0){3}{\circle*{0.15}}
\multiput(3.6,1.2)(0,0){1}{\circle*{0.15}}
\multiput(1,0.5)(4.5,0){2}{\line(1,0){0.7}}
\multiput(2.9,0.5)(0,0){1}{\line(1,0){1.4}}
\multiput(3.6,0.5)(0,0){1}{\line(0,1){0.7}}
\multiput(0,0)(6.2,0){2}{\line(1,0){1}}
\multiput(0,0)(6.2,0){2}{\line(0,1){1}}
\multiput(0,1)(6.2,0){2}{\line(1,0){1}}
\multiput(1,0)(6.2,0){2}{\line(0,1){1}}
\dashline{0.15}(1.7,0.5)(2.9,0.5)\dashline{0.15}(4.3,0.5)(5.5,0.5)
\put(0.1,0.4){$H_1$}\put(6.7,0.4){$H_2$}
\put(0.7,0.3){$v_1$}\put(6.3,0.3){$v_2$}
\put(2.3,-0.3){$i$}\put(4.9,-0.3){$j$}
\put(1.53,.42){$\underbrace{\quad\quad\quad\quad}$}
\put(4.13,.42){$\underbrace{\quad\quad\quad\quad}$}
\end{picture}
\end{center}
\caption{The graph $G_{i,j}$}\label{p6}
\end{figure}

\begin{lm}\label{lm2.8}
  Let $G_{i,j}$ be the graph shown in Figure \ref{p6} where $i,j$ are
  the numbers of included vertices.
  Then
$$\phi_{G_{i,j}}-\phi_{G_{i+1,j-1}}=(x_1-x_2)\left (p_{(H_1,v_1)}q_{(H_2,v_2)}x_2^{j-i-1}-q_{(H_1,v_1)}p_{(H_2,v_2)}x_1^{j-i-1}\right ).$$
\end{lm}

\noindent{\bf Proof }
By lemma  \ref{lm2.1}, we have
\begin{eqnarray*}
\phi_{G_{i,j}}&=&\lambda \phi_{(H_1,v_1)\cdot P_{i+j+1}\cdot (H_2,v_2)}-\phi_{(H_1,v_1)\cdot P_i}\phi_{(H_2,v_2)\cdot P_j},\\
\phi_{G_{i+1,j-1}}&=&\lambda \phi_{(H_1,v_1)\cdot P_{i+j+1}\cdot (H_2,v_2)}-\phi_{(H_1,v_1)\cdot P_{i+1}}\phi_{(H_2,v_2)\cdot P_{j-1}}.
\end{eqnarray*}
Thus, we get
\begin{eqnarray*}
\phi_{G_{i,j}}-\phi_{G_{i+1,j-1}}&=& \phi_{(H_1,v_1)\cdot P_{i+1}}\phi_{(H_2,v_2)\cdot P_{j-1}}-\phi_{(H_1,v_1)\cdot P_i}\phi_{(H_2,v_2)\cdot P_j}\\
&=&\left (p_{(H_1,v_1)}x_1^{i+1}+q_{(H_1,v_1)}x_2^{i+1}\right )\left (p_{(H_2,v_2)}x_1^{j-1}+q_{(H_2,v_2)}x_2^{j-1}\right )\\  &&-\left (p_{(H_1,v_1)}x_1^i+q_{(H_1,v_1)}x_2^i\right )\left (p_{(H_2,v_2)}x_1^j+q_{(H_2,v_2)}x_2^j\right )\\
&=& p_{(H_1,v_1)}q_{(H_2,v_2)}\left (x_1^{i+1}x_2^{j-1}-x_1^i x_2^j\right )+q_{(H_1,v_1)}p_{(H_2,v_2)}\left (x_1^{j-1}x_2^{i+1}-x_1^j x_2^i\right )\\
&=& x_1^i x_2^i \left [p_{(H_1,v_1)}q_{(H_2,v_2)}(x_1x_2^{j-i-1}-x_2^{j-i})+q_{(H_1,v_1)}p_{(H_2,v_2)}(x_1^{j-i-1}x_2-x_1^{j-i})\right ]\\
&=& (x_1-x_2)\left (p_{(H_1,v_1)}q_{(H_2,v_2)}x_2^{j-i-1}-q_{(H_1,v_1)}p_{(H_2,v_2)}x_1^{j-i-1}\right ).
\end{eqnarray*}
The proof is completed.
$\hfill\Box$

\begin{lm}\label{lm3tree}
Suppose $G_1$ and $G_2$ are two connected graphs satisfying
$G_1-u_1=G_2-u_2$ for some vertices $u_1\in V(G_1)$ and
 $u_2\in V(G_2)$. If $\phi_{G_2}(\rho(G_1))>0$, then
$\rho(G_1)>\rho(G_2)$.
\end{lm}
{\bf Proof } Let $G=G_1-u_1=G_2-u_2$. By Corollary \ref{cor0},
we have
$$ \rho(G_i) > \rho(G)
\geq \lambda_2(G_i) \quad \mbox{ for } i=1,2.$$
Here $\lambda_2(G_i)$ is the second largest eigenvalue of $G_i$. We have
$\rho(G_1)>\lambda_2(G_2)$.

Since $\rho(G_2)$ is a simple root and
$\lim\limits_{\lambda\to\infty}\phi_{G_2}(\lambda)=+\infty$, we have
 $$\phi_{G_2}(\lambda)<0 \quad \mbox{ for } \lambda\in (\lambda_2(G_2),\rho(G_2)).$$
Since $\phi_{G_2}(\rho(G_1))>0$ and $\rho(G_1)>\lambda_2(G_2)$, we must have
$\rho(G_1)>\rho(G_2)$.
\hfill $\square$

\subsection{A special tree $T_{(k-1,k,\ldots,k,k-1)}$}
The tree $T_{(k-1,k,\ldots,k,k-1)}$ ($\in \P_{n,e}$) plays an important role in this paper.
We have the following lemma.

\begin{lm}\label{lm2.9}
The spectral radius of the tree $T_{(k-1,k,...,k,k-1)}$ is the unique root $\rho_k$
of the equation $d_2=\frac{2x_1^k}{1-x_1^{k+1}}$ in the interval
$\left (\sqrt{2 +\sqrt{5}}, \infty\right )$.
\end{lm}
{\bf Remark 1:}
The following equations are equivalent to one another.
\begin{eqnarray*}
&&d_2=\frac{2x_1^k}{1-x_1^{k+1}},\\
&&d_2x_2^k-d_1x_1^k=2,\\
&&d_2=d_1x_1^{k-1},\\
&&d_2x_2^{\frac{k-1}{2}}=d_1x_1^{\frac{k-1}{2}},\\
&&d_2=2x_1^k+d_1x_1^{2k}.
\end{eqnarray*}
If ``$=$'' is replaced by ``$\geq$'', then these inequalities are still equivalent to
each other. These equivalences can be proved by equation (\ref{eq:d1d2}).
The details are omitted.

\medskip
\noindent{\bf Remark 2:} For any $k\geq 4$, we have $\rho_k\le\rho_4<\frac{3}{2}\sqrt{2}$.
For any $e\geq 6$ and $n\geq (k+2)(e-4)+6$, we can obtain a tree $T$ on $n$ vertices
and diameter $n-e$ by subdividing some edges on internal paths of $T_{(k-1,k,\ldots, k, k-1)}$.
By Lemma \ref{lm2.4}, we have
$$\rho(T)\leq \rho(T_{(k-1,k,\ldots, k, k-1)})=\rho_k< \frac{3}{2}\sqrt{2}.$$
In particular, for $e\geq 6$ and $n\geq (k+2)(e-4)+6=|T_{(k-1,k,\ldots, k, k-1)}|$,
 we have $\rho(G^{min}_{n,n-e})< \frac{3}{2}\sqrt{2}$.
In the set of graphs with spectral radius at most $\sqrt{2 +\sqrt{5}}$ (see \cite{BN}),
there is no graph with diameter $n-e$ for $e\geq 6$. Thus, $\rho(G^{min}_{n,n-e})\geq \sqrt{2+\sqrt 5}$.

\medskip
\noindent{\bf Proof of Lemma \ref{lm2.9}}.
Let $G=T_{(k-1,k,...,k,k-1)}$ and $v$ be the leftmost vertex. Note that $G$ can be built up from
a single vertex with a series of three operations as specified in Lemma
\ref{lm2.6}. We have
\begin{eqnarray*}
\phi_G&=&(1, 1) \left(\begin{array}[c]{c} p_{(T_{(k-1,k,...,k,k-1)},v)}\\ q_{(T_{(k-1,k,...,k,k-1)},v)}\end{array}\right)\\
&=& (1, 1)A^2CA^{k-1}BA^k...BA^{k-1}CA\left( \begin{array}{ll} 1 & 1\\ x_2 & x_1 \end{array}\right)^{-1}\left(\begin{array}{ll} \lambda\\ 1 \end{array}\right)\\
&=& \frac{(\lambda^2-1)^2}{x_2-x_1} (-d_1,d_2)A^{k-1}BA^k...BA^{k-1}
\left( \begin{array}{ll} d_1x_1\\ d_2x_2 \end{array}\right)\\
&=& \frac{(\lambda^2-1)^2}{x_2-x_1} (-d_1,d_2)A^{k-1}BA^k...BA^{k}
\left( \begin{array}{ll} d_1\\ d_2 \end{array}\right).
\end{eqnarray*}

Let $l=\frac{k-1}{2}$; $l$ does not have to be an integer.
Define $A^l=\left( \begin{array}{ll} x_1^l & 0\\ 0 & x_2^l \end{array}\right)$.
We can write $\phi_G$ as
\begin{equation}
  \label{eq:phiG}
\phi_G= \frac{(\lambda^2-1)^2}{x_2-x_1} (-d_1x_1^l,d_2x_2^l)(A^lBA^{l+1})^{r-1}
\left( \begin{array}{ll} d_1x_1^l\\ d_2 x_2^l\end{array}\right).
\end{equation}
It is easy to calculate
\begin{equation}\label{eq4}
A^lBA^{l+1}=\frac{1}{x_2-x_1}\left( \begin{array}{ll} d_1x_1^k & 1\\ -1 & d_2x_2^k \end{array}\right).
\end{equation}
Now we prove that $\rho_k$ is a root of $\phi_G$.
At $\lambda=\rho_k$, we have
$d_1x_1^l= d_2x_2^l$ and $d_1x_1^k+1=d_2x_2^k-1$. Thus
\begin{eqnarray*}
(A^lBA^{l+1})\left( \begin{array}{ll} 1\\ 1\end{array}\right)
& =&\frac{1}{x_2-x_1} \left( \begin{array}{ll} d_1x_1^k & 1\\ -1 & d_2x_2^k \end{array}\right)\left( \begin{array}{ll} 1\\ 1\end{array}\right)\\
&=&\frac{d_1x_1^k+1}{x_2-x_1} \left( \begin{array}{ll} 1\\ 1\end{array}\right).
\end{eqnarray*}
 We have
\begin{eqnarray*}
 \phi_G(\rho_k) &=&   \frac{(\lambda^2-1)^2}{x_2-x_1}(-d_1x_1^l, d_2x_2^l)(A^lBA^{l+1})^{r-1}\left( \begin{array}{ll} d_1x_1^l\\ d_2x_2^l \end{array}\right)\\
&=&  \frac{(\lambda^2-1)^2}{(x_2-x_1)^r}(d_1x_1^k+1)^{r-1} d_1^2x_1^{k-1} (-1,1)
\left( \begin{array}{ll} 1\\ 1\end{array}\right)\\
&=& 0.
\end{eqnarray*}

It remains to prove $\phi_G(\lambda)>0$ for any $\lambda>\rho_k$.
When $\lambda>\rho_k$, we have $d_2x_2^k-1>d_1x_1^k+1$ (and $d_2x_2^l>d_1x_1^l$).
It is easy to check
$A^lBA^{l+1}$ maps the region  $\{(z_1,z_2)\colon z_2\geq z_1>0\}$
to  $\{(z_1,z_2)\colon z_2> z_1>0\}$. By induction on $r$,
$(A^lBA^{l+1})^{r-1}$ maps the region  $\{(z_1,z_2)\colon z_2\geq z_1>0\}$
to  $\{(z_1,z_2)\colon z_2> z_1>0\}$.
Let $$\left( \begin{array}{ll} z_1\\ z_2\end{array}\right)
=(A^lBA^{l+1})^{r-1}\left( \begin{array}{ll} d_1x_1^l\\ d_2x_2^l \end{array}\right).$$
Since $d_2x_2^l>d_1x_1^l>0$, we have $z_2>z_1>0$. From equation (\ref{eq:phiG}),
we get
\begin{eqnarray*}
  \phi_G&=&\frac{(\lambda^2-1)^2}{x_2-x_1} (-d_1x_1^l,d_2x_2^l)(A^lBA^{l+1})^{r-1}
\left( \begin{array}{ll} d_1x_1^l\\ d_2 x_2^l\end{array}\right)\\
&=& \frac{(\lambda^2-1)^2}{x_2-x_1}(-d_1x_1^l,d_2x_2^l)\left( \begin{array}{ll} z_1\\ z_2\end{array}\right)\\
&=&  \frac{(\lambda^2-1)^2}{x_2-x_1}(d_2x_2^lz_2-d_1x_1^lz_1)\\
&>&0.
\end{eqnarray*}

The proof of the Lemma is finished.
$\hfill\Box$

\subsection{Limit points of some graphs}
Using the tools developed in the previous section, we can compute the limit
point of the spectral radius of some graphs.

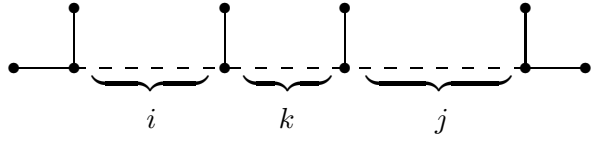
\begin{figure}[h]
\begin{center}
\setlength{\unitlength}{.8cm}
\begin{picture}(9.5,2)
\multiput(0,1)(9.5,0){2}{\circle*{0.18}}
\multiput(1,1)(0,1){2}{\circle*{0.18}}
\multiput(3.5,1)(0,1){2}{\circle*{0.18}}
\multiput(5.5,1)(0,1){2}{\circle*{0.18}}
\multiput(8.5,1)(0,1){2}{\circle*{0.18}}
\multiput(0,1)(8.5,0){2}{\line(1,0){1}}
\multiput(1,1)(2.5,0){2}{\line(0,1){1}}
\multiput(5.5,1)(3,0){2}{\line(0,1){1}}
\dashline{0.2}(1,1)(8.5,1)
\put(2.2,0){$i$}\put(4.4,0){$k$}\put(7,0){$j$}
\put(1.3,0.9){$\underbrace{\quad\quad\quad\quad}$}
\put(3.82,0.9){$\underbrace{\quad\quad\quad}$}
\put(5.85,0.9){$\underbrace{\quad\quad\quad\quad\quad}$}
\end{picture}
\end{center}
\caption{The graphs $T''_{(i,k,j)}$}\label{p7}
\end{figure}

\begin{lm}\label{lm2.10} Let $T''_{(i,k,j)}$ be the tree shown in Figure \ref{p7}
and $\rho''_{k}$ be the unique root of $d_2=x_1^{k}$ in
the interval $\left (\sqrt{2+\sqrt 5}, +\infty\right )$. Then $\lim\limits _{i,j \to \infty}\rho(T''_{(i,k,j)})=\rho_k''$.
\end{lm}
{\bf Proof } By Lemma \ref{lm2.4}, we have
$$\rho( T''_{(i,k,i)})\geq \rho(T''_{(i,k,j)}\geq T''_{(j,k,j)} \quad \mbox{ if } i\leq j.$$
It suffices to show $\lim\limits_{l\to\infty}\rho(T''_{(l,k,l)})=\rho_k''$. Let $v$ be the leftmost vertex of $T''_{(l,k,l)}$.
A simple calculation shows
\begin{eqnarray*}
&&  \phi_{T''_{(l,k,l)}}\\
&=& (1, 1) \left(\begin{array}[c]{c} p_{(T''_{(l,k,l)},v)}\\ q_{(T''_{(l,k,l)},v)}\end{array}\right)\\
 &=& (1, 1)ABA^lBA^{k}BA^lB\left( \begin{array}{ll} 1 & 1\\ x_2 & x_1 \end{array}\right)^{-1}\left( \begin{array}{ll} \lambda\\ 1 \end{array}\right)\\
&=&\frac{x_2^{2l-k+1}(d_2x_2+x_1^2)^2}{(x_2-x_1)^5}\left [((d_2x_2^{k})^2-1)-2x_1^{2l-k+3}(d_1x_1^{k}+d_2x_2^{k})-x_1^{2(2l-k+3)}((d_1x_1^{k})^2-1)\right ].
\end{eqnarray*}

As $l$ goes to infinity,   $\lim\limits_{l\to\infty}\rho(T''_{(l,k,l)})$
is the  largest root of $(d_2x_2^{k})^2-1=0$;
namely $d_2=x_1^{k}$. The proof is completed. $\hfill\Box$

We have the following Corollary  from Lemma~ \ref{lm2.10}.

\begin{cor}\label{cor2.1}
Let $T''_{(k,i)}$ be the tree shown in Figure \ref{ptki}. We have
$\lim\limits_{i \to \infty}\rho(T''_{(k,i)})=\rho''_{2k+3}.$
\end{cor}
{\bf Proof } By Lemma \ref{lm2.3}, we have
$\rho(T''_{(k,i)})=\rho(T''_{(i,2k+3,i)}).$
Thus $\lim\limits_{i \to \infty}\rho(T''_{(k,i)})=\lim\limits_{i \to \infty}
\rho(T''_{(i,2k+3,i)}) =\rho''_{2k+3}.$ \hfill $\square$

\begin{figure}[h]
\begin{center}
\setlength{\unitlength}{1.0cm}
\begin{picture}(7,2.4)
\multiput(0,.5)(0.7,0){2}{\circle*{0.15}}
\multiput(5.6,.5)(0.7,0){2}{\circle*{0.15}}
\multiput(5.6,1.2)(0,0.7){1}{\circle*{0.15}}
\multiput(.7,1.2)(0,0){1}{\circle*{0.15}}
\multiput(3.5,.5)(0,0.7){2}{\circle*{0.15}}
\multiput(.7,.5)(2.8,0){2}{\line(0,1){.7}}
\multiput(0,.5)(0,0){1}{\line(1,0){.7}}
\multiput(.7,.5)(2.8,0){2}{\line(0,1){.7}}
\multiput(5.6,.5)(0,0){1}{\line(0,1){0.7}}
\multiput(5.6,.5)(0,0){1}{\line(1,0){0.7}}
\dashline{0.15}(.7,.5)(5.6,.5)
\put(2.05,-.3){$i$}\put(4.4,-.3){$k$}
\put(0.95,.42){$\underbrace{\quad\quad\quad\quad\quad\quad}$}
\put(3.81,.42){$\underbrace{\quad\quad\quad\quad}$}
\end{picture}
\hfil
\setlength{\unitlength}{1.0cm}
\begin{picture}(7,2.4)
\multiput(0,.5)(0.7,0){2}{\circle*{0.15}}
\multiput(5.6,.5)(0.7,0){3}{\circle*{0.15}}
\multiput(5.6,1.2)(0,0.7){2}{\circle*{0.15}}
\multiput(.7,1.2)(0,0){1}{\circle*{0.15}}
\multiput(3.5,.5)(0,0.7){2}{\circle*{0.15}}
\multiput(.7,.5)(2.8,0){2}{\line(0,1){.7}}
\multiput(0,.5)(0,0){1}{\line(1,0){.7}}
\multiput(.7,.5)(2.8,0){2}{\line(0,1){.7}}
\multiput(5.6,.5)(0,0){1}{\line(0,1){1.4}}
\multiput(5.6,.5)(0,0){1}{\line(1,0){1.4}}
\dashline{0.15}(.7,.5)(5.6,.5)
\put(2.05,-.3){$j$}\put(4.4,-.3){$k$}
\put(0.95,.42){$\underbrace{\quad\quad\quad\quad\quad\quad}$}
\put(3.81,.42){$\underbrace{\quad\quad\quad\quad}$}
\end{picture}
\end{center}

\begin{multicols}{2}
\caption{The graph $T''(k,i)$}\label{ptki}\newpage
\caption{The graph $T'(k,j)$}\label{p8}
\end{multicols}
\end{figure}

\begin{lm}\label{lm2.11}
Let $T'_{(k,j)}$ be the tree shown in Figure \ref{p8} and
$\rho'_{k}$ be the unique root of $d_2=d_1^{\frac{1}{2}}x_1^{k+\frac{1}{2}}$ in
the interval $\left (\sqrt{2+\sqrt 5}, +\infty\right )$. Then $\lim\limits_{j \to \infty}\rho(T'_{(k,j)})= \rho'_{k}$.
\end{lm}
\noindent{\bf Proof } Similarly,  we have
\begin{eqnarray*}
\phi_{T'_{(k,j)}}
&=& (1, 1)  \left(\begin{array}[c]{c} p_{(T'_{(k,j)},v)}\\ q_{(T'_{(k,j)},v)}\end{array}\right)\\
&=& (1, 1)ABA^jBA^{k}CA\left( \begin{array}{ll} 1 & 1\\ x_2 & x_1 \end{array}\right)^{-1}\left( \begin{array}{ll} \lambda\\ 1 \end{array}\right)\\
&=& \frac{x_2^{j+k+1}(\lambda^2-1)(d_2x_2+x_1^3)}{(x_2-x_1)^3}\left (d_2^2-d_1x_1^{2k+1}-d_2x_1^{2j+3}-d_1^2x_1^{2j+2k+4}\right ).
\end{eqnarray*}

As $j$ goes to infinity, 
 $\lim\limits_{\j\to\infty}\rho(T'_{(k,j)})$
 is the largest root of   $d_2^2=d_1x_1^{2k+1}$;
 namely $d_2=d_1^{\frac{1}{2}}x_1^{k+\frac{1}{2}}$. The proof is  completed.
$\hfill\Box$

\subsection{Comparison of $\rho_k$,  $\rho'_{k}$, and $\rho''_{k}$}
Observe that $\rho_k$, $\rho'_k$, and $\rho''_k$ satisfy similar equations.
Since $1<\sqrt{d_1 x_1}< \frac{2}{1-x_1^{k+1}}$, we have
$$\rho_k''\leq \rho_k' \leq \rho_k.$$

For $\lambda\in [\lambda_0, \frac{3}{2}\sqrt{2}]$, $x_2$, $d_2$, and $d_1x_1$ are increasing
while $x_1$ is  decreasing. Using these facts, it is easy to check that for $k\geq 7$,
$\rho_k$, $\rho'_k$, and $\rho''_k$ are in the interval 
$(\lambda_0, \frac{3}{2}\sqrt{2})$.

We have the following lemma.
\begin{lm}\label{lm3.1}
For $k\geq 7$, we have $\rho_k<\rho''_{k-4}$ and $\rho_k<\rho'_{k-3}$.
\end{lm}
\noindent{\bf Proof } Recall that
$\rho''_{k-4}$ is the root of $d_2=x_1^{k-4}$ and $\rho_k$ is the root of $d_2=\frac{2x_1^k}{1-x_1^{k+1}}$.
We need to show  $2<x_2^4(1-x_1^{k+1})$ for $\lambda\in [\lambda_0, \frac{3}{2}\sqrt{2}]$.
For $k\geq 7$, we have
\begin{eqnarray*}
  x_2^4(1-x_1^{k+1}) &\geq& x_2^4-x_1^4\\
  &\geq&  (x_2^4-x_1^4)|_{\lambda_0} \\
  &>& 2.
\end{eqnarray*}

Note that $\rho'_{k-3}$ is the root of $d_2=\sqrt{d_1x_1}x_1^{k-3}$.
It suffices to show $2<\sqrt{d_1x_1}x_2^3(1-x_1^{k+1})$ for $\lambda\in [\lambda_0, \frac{3}{2}\sqrt{2}]$.
We have
\begin{eqnarray*}
 \sqrt{d_1x_1}x_2^3(1-x_1^{k+1}) &\geq&  \sqrt{d_1x_1}x_2^3(1-x_1^8)\\
 &\geq&  \sqrt{d_1x_1}x_2^3(1-x_1^8)|_{\lambda_0} \\
 &>& 2.
\end{eqnarray*}



The proof is completed.
$\hfill\Box$

\section{Proof of Theorem \ref{thmMain}}
The proof of Theorem \ref{thmMain} can be naturally divided into two parts.
In the first part, we prove that $G^{min}_{n,n-e}\in \P_{n,e}$.
In the second part, we prove the other statements in Theorem \ref{thmMain}.

\subsection{Part 1}


Let $\rho_{n,n-e}^{min}=\rho(G_{n,n-e}^{min})$ in the rest part of this paper.
Now we prove the following theorem, which implies the first part of Theorem \ref{thmMain}.

\begin{thm}\label{thm3.1} If $e\geq 6$ and $n\geq 10e^2-74e+142$, then
$G_{n,n-e}^{min}\in \P_{n,e}$.
 \end{thm}

\noindent{\bf Proof } By Theorem 5.2 of \cite{CDK} (see page \pageref{page2}),
it suffices to show $G_{n,n-e}^{min}\notin \P_{n,e}'$ and $G_{n,n-e}^{min}\notin \P_{n,e}''$.

Suppose $G_{n,n-e}^{min}=T'_{(k_1,k_2,...,k_{e-3})}\in \P'_{n,e}$.
Note that $T'_{(k_1,k_2,...,k_{e-3})}$ contains sub-trees of type
$T'_{(k_1,*)}$, $T''_{(k_{e-3},*)}$,  and $T''_{(*,k_i,*)}$  for $2\leq i\leq e-4$.
  By Lemma \ref{lm2.4}, Lemma \ref{lm2.10}, Corollary \ref{cor2.1},
and Lemma \ref{lm2.11},
we have
\begin{eqnarray*}
  \rho_{n,n-e}^{min} &>& \rho'_{k_1},\\
  \rho_{n,n-e}^{min} &>& \rho''_{2k_{e-3}+3},\\
\rho_{n,n-e}^{min} &>& \rho''_{k_i}, \quad \mbox{ for } 2\leq i \leq e-4.
\end{eqnarray*}
Next, we  show that at least one of $k_1, k_2,\ldots, k_{e-3}$ is small.
Let $l_1=\lceil\frac{n-3e+5}{e-3.5}\rceil$.  We claim
$$ k_1\le{l_1+1} \quad \mbox{or} \quad k_{e-3}\le\frac{l_1-3}{2}\quad \mbox{or}
\quad \exists i \in \{2,3,\ldots, e-4\}\; s.t.\; k_i \leq l_1.$$
Otherwise, we have $$k_1\geq l_1+2 \quad \mbox{and} \quad k_{e-3}\geq \frac{l_1-2}{2} \quad \mbox{and}
\quad  k_2,...,k_{e-4}
\geq l_1+1.$$
We get
$$n=\sum_{i=1}^{e-3} k_i+2e\geq l_1+2+\frac{l_1-2}{2}+(l_1+1)(e-5)+2e
=(e-3.5)l_1+3e-4
 \ge n+1.$$
Contradiction!

If $k_1\leq l_1+1$, then we have $\rho_{n,n-e}^{min} > \rho'_{l_1+1}> \rho_{l_1+4}$;
if $k_{e-3}\leq \frac{l_1-3}{2}$, then we have $\rho_{n,n-e}^{min} > \rho''_{2k_{e-3}+3}>\rho''_{l_1}> \rho_{l_1+4}$;
if $k_i\leq l_1$ for some $i\in \{2,\ldots, e-4\}$, then we have
$\rho_{n,n-e}^{min} > \rho''_{k_i}\geq \rho''_{l_1}> \rho_{l_1+4}$. In all cases,
we have
$$\rho_{n,n-e}^{min} >\rho_{l_1+4}.$$


Let $k=\lfloor \frac{n-2e+2}{e-4} \rfloor$. There exists a tree $T\in\P_{n,e}$, which
can be obtained by subdividing some edges on internal paths of $T_{(k-1,k,...,k,k-1)}$.
Since $n\geq  10e^2-74e+142$, we have
$$l_1+4=\left\lceil\frac{n-3e+5}{e-3.5}\right\rceil +4 \leq \left\lfloor
\frac{n-2e+2}{e-4}\right\rfloor=k.$$
We get
$$\rho_{n,n-e}^{min}> \rho_{l_1+4} \geq \rho(T_{(k-1,k,...,k,k-1)})\geq \rho(T).$$
Contradiction!


Now we assume $G_{n,n-e}^{min}=T''_{(k_1,k_2,...,k_{e-2})}\in \P''_{n,e}$. This is very similar
to previous case. We must have
$$k_1\le{\frac{l_2-3}{2}} \quad \mbox{or} \quad k_{e-2}\le \frac{l_2-3}{2} \quad \mbox{or} \quad \exists i\in \{2, \ldots, e-3\}\; s.t.\; k_i\le {l_2},$$
where $l_2=\lceil\frac{n-3e+7}{e-3}\rceil$.
A similar argument shows  $\rho_{n,n-e}^{min}>\rho_{l_2+4}$. Here we omit the detail.



Let $k=\lfloor \frac{n-2e+2}{e-4} \rfloor$.
There exists a tree $T\in\P_{n,e}$, which
can be obtained by subdividing some edges on internal paths of $T_{(k-1,k,...,k,k-1)}$.

Since $e\geq 5$ and $n\geq 10e^2-74e+142$, we have
$n>5e^2-31e+50$; thus,
$$l_2+4=\left\lceil\frac{n-3e+7}{e-3}\right\rceil +4 \leq \left\lfloor
\frac{n-2e+2}{e-4}\right\rfloor=k.$$
We get
$$\rho_{n,n-e}^{min}> \rho_{l_2+4} \geq \rho(T_{(k-1,k,\ldots, k,k-1)})\geq \rho(T).$$
Contradiction!

$\hfill\Box$



\noindent
{\bf Remark 3:} Assume $G^{min}_{n,n-e}=T_{(k_1,\ldots,k_r)} \in P_{n,e}$.  Let 
$\bar k=\frac{\sum_{i=1}^rk_i}{r}$. 
By Lemma \ref{lm3.1}, we can get
$k_i\ge \lfloor\overline{k}+\frac{2}{r}\rfloor-3$ for $2\le i\le r-1$ and $k_i\ge \lfloor\overline{k}+\frac{2}{r}\rfloor-2$ for $i=1,r$ whenever $n\ge 9e-30$.


\subsection{Part 2}
\begin{figure}[h]
\begin{center}
\setlength{\unitlength}{1.0cm}
\begin{picture}(7,2.4)
\multiput(0,0)(.8,0){3}{\circle*{0.15}}
\multiput(1.6,.8)(0,.8){2}{\circle*{0.15}}
\multiput(3.2,0)(0,.8){2}{\circle*{0.15}}
\multiput(5.2,0)(0,.8){2}{\circle*{0.15}}
\multiput(6.8,0)(0,.8){2}{\circle*{0.15}}
\multiput(0,0)(0,0){1}{\line(1,0){1.6}}
\multiput(1.6,0)(0,0){1}{\line(0,1){1.6}}
\multiput(3.2,0)(0,0){1}{\line(0,1){.8}}
\multiput(5.2,0)(1.6,0){2}{\line(0,1){.8}}
\dashline{0.2}(1.6,0)(6.8,0)
\multiput(3.9,.3)(0.3,0){3}{\circle*{0.1}}
\multiput(3.9,-0.3)(0.3,0){3}{\circle*{0.1}}
\put(2.25,-0.7){$k_1$}\put(5.95,-0.7){$k_j$}\put(1.65,0.1){$v_0$}\put(3.3,0.1){$v_1$}\put(5.3,0.1){$v_{j-1}$}\put(6.9,0.1){$v_j$}
\put(1.83,-0.06){$\underbrace{\quad\quad\quad}$}
\put(5.43,-0.06){$\underbrace{\quad\quad\quad}$}
\end{picture}
\end{center}
\caption{The graphs $H_{(k_1,...,k_j)}$}\label{p9}
\end{figure}
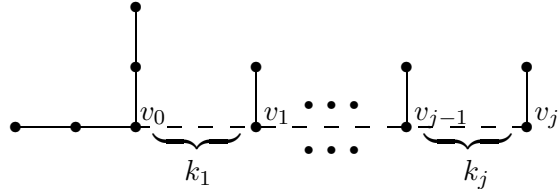

From now on, we only consider  a tree  $T_{(k_1,k_2,...,k_r)}$ in $\P_{n,e}$.
 (Here $r=e-4$ through the remaining of the paper.)
Let $v_0,v_1,...,v_r$ be the list (from left to right) of all degree $3$
vertices in $T_{(k_1,k_2,...,k_r)} \in \P_{n,e}$.
Let $H_{(k_1,k_2,...,k_j)}$ be the graph shown in Figure \ref{p9}.

Now we define two families of sub-trees of $T_{(k_1,k_2,...,k_r)}$.
For $i=1,...,r-1$, let $L_i=H_{(k_1,k_2,\ldots,k_i)}$ (from the left direction).
For $j=2,...,r$, let
 $R_j=H_{(k_r,k_{r-1},\ldots,k_j)}$ (from the right direction). We also define $L_0=P_5$ and $R_{r+1}=P_5$.

\begin{lm}\label{lm3.2}
For any $\lambda\ge{\rho(T_{(k_1,k_2,...,k_r)})}$, we have
  \begin{enumerate}
  \item $p_{(L_i,v_i)}(\lambda)\geq 0$ and $q_{(L_i,v_i)}(\lambda)\geq 0$ for $i=0,1,2,\ldots, r-1$;
  \item  $p_{(R_j,v_{j-1})}(\lambda)\geq 0$ and $q_{(R_j,v_{j-1})}(\lambda)\geq 0$ for $j=2,\ldots, r+1$.
  \end{enumerate}
\end{lm}

\noindent{\bf Proof }
For simplicity, we also write
$p_i=p_{(L_i, v_i)}$, $q_i=q_{(L_i, v_i)}$ for $i=0,1,2, \ldots, r-1$,
and $p'_j=p_{(R_j, v_{j-1})}$, $q_j'=q_{(R_j, v_{j-1})}$  for $j=2,\ldots, r+1$.
From equation (\ref{eq:P5d}),
we have $p'_{r+1}=p_0=p_{(P_5, v_0)}= \frac{d_1x_1(\lambda^2-1)}{x_2-x_1}>0$
and $q'_{r+1}=q_0=q_{(P_5, v_0)}= \frac{d_2x_2(\lambda^2-1)}{x_2-x_1}>0$
for any $\lambda> \lambda_0$.

It remains to consider
$p_i$, $q_i$ for $i=1,2,\ldots, r-1$, and $p'_j$, $q'_j$ for $j=2,\ldots, r$.
Let $\mu$ be the least number such that these functions $p_i(\lambda)$, $q_i(\lambda)$
$p'_j(\lambda)$,  $q_j'(\lambda)$ take non-negative values for all $\lambda\geq \mu$.

We need to show such $\mu$ exists.
By Lemma \ref{lmq}, we have $\lim\limits_{\lambda \to +\infty}q_i(\lambda)=+\infty$ and $\lim\limits_{\lambda \to +\infty}q_j'(\lambda)=+\infty$.
Since
$\lim\limits_{\lambda \to +\infty}p_0=\lim\limits_{\lambda \to +\infty}\frac{d_1x_1(\lambda^2-1)}{x_2-x_1}=+\infty$ and
$p_i=\frac{1}{x_2-x_1}(d_1x_1^{k_i}p_{i-1}+x_2^{k_i-1}q_{i-1})$ (see Lemma \ref{lm2.6}),
by induction on $i$, we have $\lim\limits_{\lambda \to +\infty}p_i(\lambda)=+\infty$. Similarly, we have $\lim\limits_{\lambda \to +\infty}p_j'(\lambda)=+\infty$. Thus $\mu$ is well-defined.

If $\mu\leq \rho(T_{(k_1,k_2,...,k_r)})$, then we are done. Otherwise, we assume
$\mu> \rho(T_{(k_1,k_2,...,k_r)})$.
Note that $\mu$ is always a root of one of those $p_i(\lambda)$,  $q_i(\lambda)$,
$p'_j(\lambda)$, $q_j'(\lambda)$.

\begin{description}

\item[Case (1)] There exists an $i$ ($1\leq i\leq r-1$) such that $p_i(\mu)=0$.
 Since  $p_i=\frac{1}{x_2-x_1}(d_1x_1^{k_i}p_{i-1}+x_2^{k_i-1}q_{i-1})$, we must have $p_{i-1}(\mu)=q_{i-1}(\mu)=0.$
By Lemma \ref{lm2.7},
we have
$$\phi_{T_{(k_1,k_2,...,k_r)}}(\mu)=(x_2-x_1)(x_2^{k_i-1}q_{i-1}q_{i+1}'-x_1^{k_i-1}p_{i-1}p_{i+1}')\mid_{\mu}=
0.$$
It contradicts to the assumption $\mu> \rho(T_{(k_1,k_2,...,k_r)})$.

\item[Case (2)] There exists a $j$ ($2\leq j\leq r$)  such that $p_j'(\mu)=0$.
This case is symmetric to Case (1).

\item[Case (3)] There exists an $i$ ($1\leq i\leq r-1$) such that $q_i(\mu)=0$.
By Lemma \ref{lm2.7},  we have
$$\phi_{T_{(k_1,k_2,...,k_r)}}(\mu)=(x_2-x_1)(x_2^{k_{i+1}-1}q_{i}q_{i+2}'-x_1^{k_{i+1}-1}p_{i}p_{i+2}')\mid_{\mu}\leq
0.$$
It contradicts to $\mu> \rho(T_{(k_1,k_2,...,k_r)})$.

\item[Case (4)] There exists a $j$ ($2\leq j\leq r$)  such that $q_j'(\mu)=0$. This case is symmetric to Case (3).
\end{description}
The proof of this Lemma is finished.
\hfill $\square$



The following Lemma gives the lower bound
for the spectral radius of a general tree $T_{(k_1,k_2,...,k_r)}\in\P_{n,e}$.
\begin{lm}\label{lm3.3} 
 Let $\overline k=\frac{\sum_{i=1}^r k_i}{r}$.
We have
$$d_2\ge\frac{2x_1^{\overline{k}+\frac{2}{r}}}{1-x_1^{\overline{k}+\frac{2}{r}+1}}$$  for all $\lambda\ge\rho(T_{(k_1,k_2,...,k_r)})$, where the equality
holds if and only if $k_1+1=k_2=\cdots=k_{r-1}=k_r+1$ and $\lambda =\rho(T_{(k_1,k_2,...,k_r)})$.
\end{lm}

\noindent{\bf Proof } For $i=0,1,2, \ldots, r-1$, we define
$t_i=q_i/p_i$. Similarly, for $j=2, \ldots, r+1$, we define
$t'_j=q'_j/p'_j$.
For any $s>0$, we define $$f_s(t)=\frac{d_2x_2^{2s}t-x_2}{x_2^{2s-1}t+d_1}=\frac{d_2x_2t-x_1^{2s-2}}{t+d_1x_1^{2s-1}},\quad t>0.$$
We consider the fixed point of $f_s(t)$, which satisfies
$$t^2-(d_2x_2-d_1x_1^{2s-1})t+x_1^{2s-2}=0.$$

This quadratic equation has a unique root $x_1^{s-1}$ when
\begin{equation}
  \label{eq:slambda}
d_2=2x_1^{s}+d_1x_1^{2s}.
\end{equation}
We choose $s=s(\lambda)$ to be the root of Equation (\ref{eq:slambda}). The line $y=t$ is tangent to the
curve $y=f_s(t)$ at $t=x_1^{s-1}$. Because $f_s(t)$ is an increasing and concave function of $t$,
we have
$$f_s(t)\le t,  \quad \forall t>0.$$

For  $i=1,...,r$, we have
\begin{equation}
  \label{eq:kit}
f_{k_i}(t)=f_{s}(x_2^{2(k_i-s)}t)\le{x_2^{2(k_i-s)}t}.
\end{equation}

By Lemma \ref{lm2.7},  we get
$$\phi_{T_{(k_1,k_2,...,k_r)}}=(x_2-x_1)(x_2^{k_r-1}q_{r-1}q'_{r+1}-x_1^{k_r-1}p_{r-1}p'_{r+1}).$$
Since $\phi_{T_{(k_1,k_2,...,k_r)}}\geq 0$ for all $\lambda\ge{\rho(T_{(k_1,k_2,...,k_r)})}$, we get
 $$t_{r-1}t'_{r+1}x_2^{2(k_r-1)}\ge 1.$$
Note $t'_{r+1}=t_0=\frac{d_2x_2}{d_1x_1}=\frac{d_2}{d_1}x_2^2$. 
Applying inequality (\ref{eq:kit}) recursively, we have
\begin{eqnarray*}
1&\leq& \frac{d_2}{d_1}x_2^2\cdot x_2^{2(k_r-1)} \frac{q_{r-1}}{p_{r-1}}\\
&=& \frac{d_2}{d_1}x_2^{2k_r}f_{k_{r-1}}(f_{k_{r-2}}(...(f_{k_1}(t_0)...)))\\
&\le &  \frac{d_2}{d_1} x_2^{2k_r}  x_2^{2(k_{r-1}-s)}x_2^{2(k_{r-2}-s)}...x_2^{2(k_1-s)}t_0 \\
&=&\frac{d_2}{d_1} x_2^{2k_r}  x_2^{2(k_{r-1}-s)}x_2^{2(k_{r-2}-s)}...x_2^{2(k_1-s)} \frac{d_2}{d_1}x_2^2\\
&=& \frac{d_2^2}{d_1^2}x_2^{2(r\overline{k}-(r-1)s+1)}.
\end{eqnarray*}
We get $d_2\geq d_1x_1^{r\overline{k}-(r-1)s+1}$; and the equality holds if and only if $k_1+1=k_2=\cdots=k_{r-1}=k_r+1=s$  and $\lambda=\rho(T_{(k_1,k_2,\ldots, k_r)})$.
By Remark 1, $d_2\geq d_1x_1^{r\overline{k}-(r-1)s+1}$ is equivalent to
\begin{equation}
  \label{eq:sd2}
d_2\geq 2x_1^{r\overline{k}-(r-1)s+2}+ d_1x_1^{2(r\overline{k}-(r-1)s+2)}.
\end{equation}
Comparing this inequality with equation (\ref{eq:slambda}), we must have
 $s\le r\overline{k}-(r-1)s+2$. Solving $s$, we get $s\le \overline{k}+\frac{2}{r}$.
Thus, $$d_2=2x_1^s+d_1x_1^{2s}\ge 2x_1^{  \overline{k}+\frac{2}{r}}+d_1x_1^{2( \overline{k}+\frac{2}{r})}.$$
Applying Remark 1 one more time, we get $$d_2\ge\frac{2x_1^{ \overline{k}+\frac{2}{r}}}{1-x_1^{\overline{k}+\frac{2}{r}+1}}.$$
The proof is completed.  $\hfill\Box$

\begin{lm}\label{lm3.4} Let $G_{n,n-e}^{min}=T_{(k_1,k_2,...,k_r)}$
  and $\overline{k}=\frac{\sum_{i=1}^{r} k_i}{r}$.
  Then $$d_2\le\frac{2x_1^{\lfloor\overline{k}+\frac{2}{r}\rfloor}}{1-x_1^{\lfloor\overline{k}+\frac{2}{r}\rfloor+1}}$$ holds at $\lambda=\rho_{n,n-e}^{min}$.
\end{lm}

\noindent{\bf Proof } Let $s=\overline{k}+\frac{2}{r}$. Observe that we can always subdivide some edges on internal paths of
$T_{(\lfloor s \rfloor-1,\lfloor s \rfloor,...,\lfloor s \rfloor,\lfloor s \rfloor-1)}$ to get a tree $T$ on $n$ vertices and diameter
$n-e$. By Lemma \ref{lm2.4},
we have
$$\rho_{n,n-e}^{min}\le \rho(T) \leq \rho(T_{(\lfloor s \rfloor-1,\lfloor s \rfloor,...,\lfloor s \rfloor,\lfloor s \rfloor-1)}) =\rho_{\lfloor s \rfloor}.$$

By Lemma  \ref{lm2.9}, $\rho_{\lfloor s \rfloor}$ is the root of $$d_2=\frac{2x_1^{\lfloor s \rfloor}}{1-x_1^{\lfloor s \rfloor+1}}.$$
Since $d_2(\lambda)$ is increasing while $\frac{2x_1^{\lfloor s \rfloor}}{1-x_1^{\lfloor s\rfloor+1}}$ is decreasing on $\left (\sqrt{2+\sqrt{5}}, \infty\right )$, we get
$$d_2(\rho_{n,n-e}^{min})\le d_2(\rho_{\lfloor s \rfloor})=
\left .\frac{2x_1^{\lfloor s \rfloor}}{1-x_1^{\lfloor s\rfloor+1}}\right |_{\rho_{\lfloor s \rfloor}} \le\left .
\frac{2x_1^{\lfloor s \rfloor}}{1-x_1^{\lfloor s\rfloor+1}}\right |_{\rho_{n,n-e}^{min}}.$$

The proof is completed. $\hfill\square$

We get the following corollary.

\begin{cor}\label{thm3.2} Let $G_{n,n-e}^{min}=T_{(k_1,k_2,...,k_r)}\in \P_{n,e}$ and $s=\frac{1}{r}\sum_{i=1}^r k_i+\frac{2}{r}= \frac{n-2e+2}{e-4}$.
We have $$\frac{2x_1^s}{1-x_1^{s+1}}\le{d_2}\le\frac{2x_1^{\lfloor s \rfloor}}{1-x_1^{\lfloor s\rfloor+1}}$$ holds at $\lambda=\rho(G_{n,n-e}^{min})$.
In particular, $\rho(G_{n,n-e}^{min})=\sqrt{2+\sqrt{5}}+ O\left((\frac{\sqrt{5}-1}{2})^{s/2}\right).$
\end{cor}



\begin{lm}\label{lm3.5}
Assume  $G_{n,n-e}^{min}=T_{(k_1,...,k_i,k_{i+1},...,k_{r})}$ and $\overline{c}=\frac{\rho_{n,n-e}^{min}+\sqrt{(\rho_{n,n-e}^{min})^2+4d_1d_2}}{2}$.
Then the following equalities hold at the point $\lambda=\rho_{n,n-e}^{min}$.
\begin{eqnarray}
  \label{eq:c12}
  \overline{c}x_1^{k_i+1}\le &d_2&\le\overline{c}x_1^{k_i-1} \quad \mbox{ for } i=2,...,r-1;\\
\label{eq:c13}
\sqrt{\overline{c}d_1}x_1^{k_i+1}\le &d_2&\le \sqrt{\overline{c}d_1}x_1^{k_i} \quad \mbox{ for } i=1,r.
\end{eqnarray}

\end{lm}
\noindent{\bf Proof } We reuse notations
 $L_i, p_i, q_i, t_i $ (for $i=0,1,...,r-1$) and
$R_j, p'_j, q'_j, t'_j$ (for $j=2,...,r+1$), which are
introduced in Lemma \ref{lm3.2}  and
 Lemma \ref{lm3.3}.

Choosing any $i\in{\{1,...,r-1\}}$, by Lemma  \ref{lm2.7} we have $$t_it'_{i+2}x_2^{2(k_{i+1}-1)}=1$$ at $\lambda=\rho_{n,n-e}^{min}$. This means
$$\frac{d_2x_2t_{i-1}-x_1^{2(k_i-1)}}{t_{i-1}+d_1x_1^{2k_i-1}}t'_{i+2}x_2^{2(k_{i+1}-1)}=1.$$
We can rewrite it as
\begin{equation}\label{eq:product}
\left (t_{i-1}-\frac{x_1^{2k_i-1}}{d_2}\right )\left (t'_{i+2}-\frac{x_1^{2k_{i+1}-1}}{d_2}\right )=\frac{d_1d_2+1}{d_2^2}x_1^{2(k_i+k_{i+1}-1)}.
\end{equation}
Note $t_i=f_{k_i}(t_{i-1})=\frac{d_2x_2t_{i-1}-x_1^{2(k_i-1)}}{t_{i-1}+d_1 x_1^{2k_i-1}}>0$.
We have
\begin{equation}
  \label{eq:ti-1}
  t_{i-1}>\frac{x_1^{2k_i-1}}{d_2}.
\end{equation}

For $i=1,\ldots, r-1$, we apply Lemma \ref{lm3tree} to
$G_1=T_{(k_1,\ldots, k_{i}, k_{i+1},\ldots, k_r)}$ and $G_2=T_{(k_1,\ldots, k_{i}+1, k_{i+1}-1,\ldots, k_r)}$,
where both trees contain a common induced subtree
$T_{(k_1,\ldots, k_{i}+k_{i+1}+1,\ldots, k_r)}$ (after removing one leaf vertex).
If $\phi_{G_2}(\rho(G_1))>0$,  then $\rho(G_1)> \rho(G_2)$.
This contradict to the assumption $G_1=G^{min}_{n,n-e}$.

We get $\phi_{G_2}(\rho(G_1))\leq 0$, i.e., $\phi_{T_{(k_1,\ldots, k_{i}+1, k_{i+1}-1,\ldots, k_r)}}
(\rho_{n,n-e}^{min})\leq 0.$

We apply Lemma  \ref{lm2.8} to obtain  the difference of characteristic polynomials of $T_{(k_1,\ldots, k_{i}, k_{i+1},\ldots, k_r)}$
and $T_{(k_1,\ldots, k_{i}+1, k_{i+1}-1,\ldots, k_r)}$,
\[
\phi_{T_{(k_1,\ldots, k_{i}, k_{i+1},\ldots, k_r)}}-\phi_{T_{(k_1,\ldots, k_{i}+1, k_{i+1}-1,\ldots, k_r)}} =
(x_1-x_2)\left (p_{i-1}q_{i+2}' x_2^{k_{i+1}-k_i-1}-q_{i-1}p_{i+2}' x_1^{k_{i+1}-k_i-1}\right ).
 \]
Evaluating the function above  at  $\lambda=\rho^{min}_{n,n-e}$, we have
$$(x_1-x_2)\left.\left (p_{i-1}q_{i+2}' x_2^{k_{i+1}-k_i-1}-q_{i-1}p_{i+2}' x_1^{k_{i+1}-k_i-1}\right )\right |_{\rho_{n,n-e}^{min}}\ge 0.$$
Since $q_{i-1}\geq 0$ and $p_{i+2}'\geq 0$ (from Lemma \ref{lm3.2}),
we get $\frac{t'_{i+2}}{t_{i-1}}\le{x_1^{2(k_{i+1}-k_i-1)}}$ at $\lambda=\rho_{n,n-e}^{min}.$
In the rest of the proof, all expressions are evaluated at $\lambda=\rho_{n,n-e}^{min}$.
The notation ``$|_{\rho_{n,n-e}^{min}}$'' is omitted for simplicity.

On the one hand,  by inequality (\ref{eq:ti-1}),
 we can substitute ${t'_{i+2}}\le{t_{i-1}}{x_1^{2(k_{i+1}-k_i-1)}}$ into equation(\ref{eq:product}) 
and get
$$\left (t_{i-1}-\frac{x_1^{2k_i-1}}{d_2}\right )\left (x_1^{2(k_{i+1}-k_i-1)}t_{i-1}-\frac{x_1^{2k_{i+1}-1}}{d_2}\right )\ge\frac{d_1d_2+1}{d_2^2}x_1^{2(k_i+k_{i+1}-1)}.$$
After simplification, we have
\[
d_2t_{i-1}^2-x_1^{2k_i}\rho_{n,n-e}^{min}t_{i-1}-d_1x_1^{4k_i}\ge 0.
\]
Recall $\bar c=\frac{\rho_{n,n-e}^{min}+\sqrt{(\rho_{n,n-e}^{min})^2+4d_1d_2}}{2}$.
Solving this quadratic inequality, since $t_{i-1}>0$, we get
$$t_{i-1}\ge \bar c x_1^{2k_i} / d_2, \quad i=1,...,r-1.$$

By symmetry, we have $$t'_{i+1}\ge \bar c x_1^{2k_i}/d_2, \quad i=2,...,r.$$

On the other hand, we substitute ${t_{i-1}}\ge{t'_{i+2}}{x_2^{2(k_{i+1}-k_i-1)}}$ into equation (\ref{eq:product}).
By the similar calculation, we get

$$t'_{i+2}\le \bar c x_1^{2(k_{i+1}-1)}/d_2, \quad i=1,...,r-1.$$
Changing the index $i+2$ to $i+1$, we have
 $$t'_{i+1} \le \bar c x_1^{2(k_i-1)}/d_2, \quad i=2,...,r. $$
By symmetry, we have
$$t_{i-1}\le \bar c x_1^{2(k_i-1)}/d_2, \quad i=1,...,r-1.$$

Combining the inequalities above, we get
\begin{eqnarray}
  \label{eq:c15}
\frac{\bar c}{d_2}x_1^{2k_i} &\le& t_{i-1}\le \frac{\bar c}{d_2}x_1^{2(k_i-1)}, \quad i=1,...,r-1,\\
\label{eq:c16}
\frac{\bar c}{d_2}x_1^{2k_i} &\le& t'_{i+1}\le \frac{\bar c}{d_2}x_1^{2(k_i-1)}, \quad i=2,...,r.
\end{eqnarray}

Now we apply Lemma  \ref{lm2.7} and get
\begin{equation}
  \label{eq:c17}
  t_{i-1}t'_{i+1}x_2^{2(k_i-1)}= 1.
\end{equation}
Taking product of inequalities (\ref{eq:c15}),  (\ref{eq:c16}), and then substituting $t_{i-1}t'_{i+1}$ into equation (\ref{eq:c17}).
After simplification, wet get inequality (\ref{eq:c12}).

When $i=1$ or $r$, we have
$$\frac{\bar c}{d_2}x_1^{2k_i}\le t_0=t'_{r+1}=\frac{d_2}{d_1}x_2^2\le \frac{\bar c}{d_2} x_1^{2(k_i-1)}. $$
Solving for $d_2$, we get inequality  (\ref{eq:c13}).
The proof of this lemma is completed. $\hfill\Box$
\\


\noindent{\bf Proof of the second part of theorem \ref{thmMain}}.
As in the proof of  Lemma \ref{lm3.5},  all expressions in this proof
are evaluated at $\lambda=\rho_{n,n-e}^{min}$
and $``|_{\rho_{n,n-e}^{min}}"$ is omitted for simplicity.

By Lemma \ref{lm3.5}, for $2\le i\le r-1$, we have
$$\overline{c}x_1^{k_i+1}\le d_2\le\overline{c}x_1^{k_i-1}.$$
By the definition of $\bar c$, we get
$$2d_2x_1^{k_i+1}\le\rho_{n,n-e}^{min}+\sqrt{(\rho_{n,n-e}^{min})^2+4d_1d_2}\le2d_2x_1^{k_i-1}.$$
After solving $d_2$ and simplifying, we have $$\frac{\rho_{n,n-e}^{min}x_1^{k_i+1}+2x_1^{2k_i+3}}{1-x_1^{2(k_i+2)}}\le{d_2}\le\frac{\rho_{n,n-e}^{min}x_1^{k_i-1}+2x_1^{2k_i-1}}{1-x_1^{2k_i}}.$$

Since $\rho_{n,n-e}^{min}>2>1+x_1^2=\rho_{n,n-e}^{min}x_1$, we observe
 $$\frac{2x_1^{k_i+1}}{1-x_1^{k_i+2}}<\frac{\rho_{n,n-e}^{min}x_1^{k_i+1}+2x_1^{2k_i+3}}{1-x_1^{2(k_i+2)}}$$
and $$\frac{\rho_{n,n-e}^{min}x_1^{k_i-1}+2x_1^{2k_i-1}}{1-x_1^{2k_i}}<\frac{2x_1^{k_i-2}}{1-x_1^{k_i-1}}.$$

We obtain
\begin{equation}
  \label{eq:c18}
  \frac{2x_1^{k_i+1}}{1-x_1^{k_i+2}}<d_2<\frac{2x_1^{k_i-2}}{1-x_1^{k_i-1}} \quad \mbox{ for }2\le i \le r-1.
\end{equation}
From Theorem \ref{thm3.2}, we have
\begin{equation}
  \label{eq:c19}
\frac{2x_1^s}{1-x_1^{s+1}}\le{d_2}\le\frac{2x_1^{\lfloor s \rfloor}}{1-x_1^{\lfloor s\rfloor+1}}.
\end{equation}
Combining inequalities (\ref{eq:c18}) and (\ref{eq:c19}), we get
\begin{eqnarray*}
\frac{2x_1^{k_i+1}}{1-x_1^{k_i+2}}&<&\frac{2x_1^{\lfloor s \rfloor}}{1-x_1^{\lfloor s\rfloor+1}},\\
\frac{2x_1^{k_i-2}}{1-x_1^{k_i-1}}&>&\frac{2x_1^s}{1-x_1^{s+1}}.
\end{eqnarray*}
Thus, $\lfloor s \rfloor-1<k_i<s+2$. So $\lfloor s \rfloor<k_i\le\lceil s \rceil+1$ where $i=2,...,r-1$.

For $j=1$ or $r$, combining inequalities (\ref{eq:c13}) and  (\ref{eq:c19}),
we have
\begin{eqnarray*}
\frac{\sqrt{\overline{c}d_1}}{2}x_1^{k_j+1}&\le&\frac{x_1^{\lfloor s \rfloor}}{1-x_1^{\lfloor s\rfloor+1}},\\
\frac{\sqrt{\overline{c}d_1}}{2}x_1^{k_j}&\ge&\frac{x_1^s}{1-x_1^{s+1}}.
\end{eqnarray*}

Note that $d_1 \to 2x_1$ and  $\overline{c} \to \lambda_0$ as $n$ approaches infinity.
For sufficiently large $n$, we have
$x_2^{0.1}<\frac{\lambda_0}{2}<x_2^{0.2}$. We get $$x_1^{k_j+1+0.45}\le x_1^{\lfloor s \rfloor} \quad \mbox{and} \quad x_1^{k_j+0.4}>x_1^s.$$
So $\lfloor s \rfloor-1\le k_j\le\lfloor s \rfloor$ for $n$ large enough.

In conclusion, we get
$$\lfloor s \rfloor-1\le k_j\le\lfloor s \rfloor \leq k_i \leq \lceil s \rceil +1$$
for $2\leq i\leq r-1$ and $j=1,r$.

Now we will prove item 2. It suffices to show $k_i-k_j\leq 2$, for $2\leq i\leq r-1$ and $j=1,r$.
Suppose that there exist $i,j$ with $i\in \{2,...,r-1\}$ and $j\in \{1,r\}$ so that $k_i\ge {k_j+3}$.
By Lemma \ref{lm3.5}, we have $$\sqrt{\overline{c}d_1}x_1^{k_j+1}\le d_2\le{\overline{c}x_1^{k_j+2}}.$$

Since $\lambda x_1^2=(1+x_1^2)x_1<2x_1\le d_1$ for $\lambda\ge\lambda_0$ and $\overline{c}\to \lambda_0$ as $n$ approaches infinity, we have
$\overline{c}x_1^{k_j+2}<\sqrt{\overline{c}d_1} x_1^{k_j+1}$ for $n$ large enough. Contradiction!


Now we will prove item 3. By Lemma \ref{lm3.5}, we have
$\overline{c}x_1^{k_j+1}\le d_2\le \overline{c}x_1^{k_i-1}$ for all $2\le i,j\le r-1$.  This implies
$|k_i-k_j|\le 2$.
It is sufficient to show that there are no $i,j$ with $|k_i-k_j|=2$. Otherwise,
suppose there exist $i,j  \in \{2,...,r-1\}$ such that $k_i=k$ and $k_j=k+2$.
Without loss of generality, we can assume that  $i<j$ and in addition $i,j$ are mostly close to each other.
Namely, $k_l=k+1$ for all integer $l$ between $i$ and $j$.

Applying inequality (\ref{eq:c12}) to $k_i=k$ and $k_j=k+2$, we have
\begin{eqnarray*}
  d_2 &\geq& \overline{c}x_1^{k_i+1} = \overline{c}x_1^{k+1},\\
  d_2 &\leq&  \overline{c}x_1^{k_j-1} = \overline{c}x_1^{k+1}.
\end{eqnarray*}
Two inequalities above force $d_2=\overline{c}x_1^{k+1}$. These equalities force
$t_{i-1}=t'_{i+1}=x_1^{k-1}$, $t_{j-1}=t'_{j+1}=x_1^{k+1}$ by inequalities (\ref{eq:c15}) and (\ref{eq:c16}).

Consider the function $f(t)=\frac{d_2t-x_2}{x_1t+d_1}={f_k(x_1^{2k}t)}$ and let $c=\bar c/d_2=x_2^{k+1}$.
It is easy to check $f(c)=\frac{1}{c}$. We claim
$$t_l=x_1^{k+1} \quad \mbox{ for } i\leq l \leq j-1.$$
For $l=i$, we have
$$t_i=f_k(t_{i-1})=f_k(x_1^{k-1})=f(x_2^{k+1})=f(c)=\frac{1}{c}=x_1^{k+1}.$$
By induction on $l$, we have
$$t_l=f_{k+1}(t_{l-1})=
f_{k+1}(x_1^{k+1})=f(x_2^{k+1})=f(c)=\frac{1}{c}=x_1^{k+1}.$$
By Lemma \ref{lm2.7}, we have
$$ t_{j-2}t_j'x_2^{2k}=1.$$
Since $t_{j-2}=x_1^{k+1}$, it implies $t_j'=x_1^{k-1}$.
However, we also have
$$t_j'=f_{k+2}(t_{j+1}')= f_{k+2}(x_1^{k+1})=f(x_2^{k+3})\not= x_1^{k-1}.$$
Contradiction!

If $n-6$ is divisible by $e-4$, then $s=\frac{n-6}{e-2}-4$ is an
integer.  In this case, the only possible sequence $(k_1, k_2,\ldots,
k_r)$ satisfying items 1-3 is $(s-1, s,\ldots, s, s-1)$. In
particular, we have $G^{min}_{n,n-e}=T_{(s-1, s,\ldots, s, s-1)}$.

The proof is completed.
$\hfill\Box$





\section{Proofs of Theorems \ref{thm1.3} and \ref{thm1.4}}
\subsection{e=7}
Let $G_{n,n-7}^{min}=T_{(k_1,k_2,k_3)}\in \P_{n,7}$. Note $k_1+k_2+k_3=n-14$.
By Theorem \ref{thmMain}, here are all the possible graphs for $G_{n,n-7}^{min}$.
\begin{description}
\item[Case 1.] $\sum\limits_{i=1}^3 k_i=3k$. We have $(k_1,k_2,k_3)=(k,k,k)$ or $(k,k+1,k-1)$.
\item[Case 2.] $\sum\limits_{i=1}^3 k_i=3k+1$. We have $(k_1,k_2,k_3)=(k,k+1,k)$.
\item[Case 3.] $\sum\limits_{i=1}^3 k_i=3k+2$. We have $(k_1,k_2,k_3)=(k,k+2,k)$ or $(k,k+1,k+1)$.
\end{description}

To simplify the proof of Theorem \ref{thm1.3}, we introduce the following short notations.
We have
\begin{eqnarray*}
  p_0&:=&p_{(L_0,v_0)}=\frac{\lambda^2-1}{x_2-x_1}d_1x_1, \\
  q_0&:=&q_{(L_0,v_0)}=\frac{\lambda^2-1}{x_2-x_1}d_2x_2, \\
 p^{(k-1)}&:=&p_{(H_{(k-1)},v_1)}=\frac{\lambda^2-1}{(x_2-x_1)^2}(d_1^2x_1^k+d_2x_2^{k-1}),\\
 q^{(k-1)}&:=& q_{(H_{(k-1)},v_1)} =\frac{\lambda^2-1}{(x_2-x_1)^2}(d_2^2x_2^k-d_1x_1^{k-1}),\\
  p^{(k)}&:=&p_{(H_{(k)},v_1)} = \frac{\lambda^2-1}{(x_2-x_1)^2}(d_1^2x_1^{k+1}+d_2x_2^k), \\
  q^{(k)}&:=&q_{(H_{(k)},v_1)} =  \frac{\lambda^2-1}{(x_2-x_1)^2}(d_2^2x_2^{k+1}-d_1x_1^k), \\
  p^{(k+1)}&:=& p_{(H_{(k+1)},v_1)}=\frac{\lambda^2-1}{(x_2-x_1)^2}(d_1^2x_1^{k+2}+d_2x_2^{k+1}),\\
  q^{(k+1)}&:=& q_{(H_{(k+1)},v_1)}=\frac{\lambda^2-1}{(x_2-x_1)^2}(d_2^2x_2^{k+2}-d_1x_1^{k+1}),\\
  p^{(k,k+1)}&:=&p_{(H_{(k,k+1)},v_2)}=\frac{\lambda^2-1}{(x_2-x_1)^3}(d_1^3x_1^{2k+2}+d_1d_2x_1+d_2^2x_2^{2k+1}-d_1),\\
 q^{(k,k+1)}&:=&q_{(H_{(k,k+1)},v_2)}=\frac{\lambda^2-1}{(x_2-x_1)^3}(d_2^3x_2^{2k+2}-d_1d_2x_2-d_1^2x_1^{2k+1}-d_2).
\end{eqnarray*}




\noindent{\bf Proof of Theorem \ref{thm1.3}}. We will compare the spectral radius of the possible graphs listed above in three cases separately.

\textbf{Case 1.} $\sum\limits_{i=1}^3 k_i=3k$.

By Lemma \ref{lm2.8}, we have
\begin{eqnarray*}
\phi_{T_{(k,k,k)}}-\phi_{T_{(k,k+1,k-1)}}&=&(x_1-x_2)\left (p^{(k)}q_0x_1-q^{(k)}p_0x_2\right )\\
&=&-\frac{(\lambda^2-1)^2}{(x_2-x_1)^2}\left [(d_2x_1+1)d_1^2x_1^k-(d_1x_2-1)d_2^2x_2^k\right ]\\
&=&\frac{(d_2x_1+1)(\lambda^2-1)^2}{(x_2-x_1)^2}\left (d_2^2x_2^k-d_1^2x_1^k\right ).
\end{eqnarray*}
In the last step, we applied the fact $d_2x_1+1=d_1x_2-1$.

By Lemma \ref{lm2.9} and Remark 1, $\rho(T_{(k,k+1,k)})$ ($=\rho_{k+1}$) satisfies $d_2x_2^{k/2}=d_1x_1^{k/2}$. The largest root of
 $\phi_{T_{(k,k,k)}}-\phi_{T_{(k,k+1,k-1)}}=0$ is $\rho_{k+1}$.


Noting that $d_2^2x_2^k-d_1^2x_1^k$ is an increasing function of $\lambda\in\left (\sqrt{2+\sqrt 5}, \frac{3}{2}\sqrt 2\right )$ for sufficiently large $k$. 
 By Lemma \ref{lm2.4}, we have
$\rho_{k+1}=\rho(T_{(k,k+1,k)})<\rho(T_{(k,k,k)})$.  Evaluating
 $\phi_{T_{(k,k,k)}}-\phi_{T_{(k,k+1,k-1)}}$
at $\lambda=\rho(T_{(k,k,k)})$, we get  $\phi_{T_{(k,k+1,k-1)}}(\rho(T_{(k,k,k)}))<0$.
Thus, by Lemma \ref{lm2.2}, $\rho(T_{(k,k,k)})<\rho(T_{(k,k+1,k-1)})$ and $G_{n,n-7}^{min}=T_{(k,k,k)}$.

\textbf{Case 2.} $\sum\limits_{i=1}^3 k_i=3k+1$.  We must have $G_{n,n-7}^{min}=T_{(k,k+1,k)}$.

\textbf{Case 3.} $\sum\limits_{i=1}^3 k_i=3k+2$.

Similarly by Lemma \ref{lm2.8}, we have $$\phi_{T_{(k,k+1,k+1)}}-\phi_{T_{(k,k+2,k)}}=\frac{(d_2x_1+1)(\lambda^2-1)^2}{(x_2-x_1)^2}\left (d_2^2x_2^k-d_1^2x_1^k\right ).$$
Noting that $d_2^2x_2^k-d_1^2x_1^k$ is an increasing function of $\lambda\in\left (\sqrt{2+\sqrt 5}, \frac{3}{2}\sqrt 2\right )$ for sufficiently large $k$. We have
$\phi_{T_{(k,k+1,k+1)}}(\lambda)<\phi_{T_{(k,k+2,k)}}(\lambda)$ for
any $\sqrt{2+\sqrt 5}\leq \lambda< \rho_{k+1}$.
By Lemma \ref{lm2.4}, we get $\rho(T_{(k,k+2,k)})<\rho(T_{(k,k+1,k)})=\rho_{k+1}$.
 Thus, $\phi_{T_{(k,k+1,k+1)}}(\rho(T_{(k,k+2,k)}))<0$.
It follows $\rho(T_{(k,k+1,k+1)})>\rho(T_{(k,k+2,k)})$. So $G_{n,n-7}^{min}=T_{(k,k+2,k)}$.

The proof of Theorem \ref{thm1.3} is completed.
$\hfill\Box$

\subsection{e=8}

Now we let $G_{n,n-8}^{min}=T_{(k_1,k_2,k_3,k_4)}\in \P_{n,8}$. By Theorem \ref{thmMain}, all the possible graphs for $G_{n,n-8}^{min}$ are as follows.
\begin{description}
\item[Case 1.] If $\sum\limits_{i=1}^4 k_i=4k$, then $(k_1,k_2,k_3,k_4)=(k,k,k,k)$, $(k,k,k+1,k-1)$,  $(k,k+1,k,k-1)$, or $(k-1,k+1,k+1,k-1)$.
\item[Case 2.] If $\sum\limits_{i=1}^4 k_i=4k+1$, then $(k_1,k_2,k_3,k_4)=(k,k+1,k,k)$ or $(k,k+1,k+1,k-1)$.
\item[Case 3.] If $\sum\limits_{i=1}^4 k_i=4k+2$, then $(k_1,k_2,k_3,k_4)=(k,k+1,k+1,k)$.
\item[Case 4.] If $\sum\limits_{i=1}^4 k_i=4k+3$, then $(k_1,k_2,k_3,k_4)=(k,k+1,k+1,k+1)$ or $(k,k+1,k+2,k)$.
\end{description}



\noindent{\bf Proof of Theorem \ref{thm1.4}}. Similarly, we denote
$p^{(k,k)}=P_{(H_{(k,k)},v_2)}$, $q^{(k,k)}=q_{(H_{(k,k)},v_2)}$,
$p^{(k-1,k+1)}=P_{(H_{(k-1,k+1)},v_2)}$, and $q^{(k-1,k+1)}=q_{(H_{(k-1,k+1)},v_2)}$.

We will compare the spectral radius of all possible graphs listed in four cases above.

\textbf{Case 1.} $\sum\limits_{i=1}^4 k_i=4k$.

First we prove
$$\rho(T_{(k,k,k,k)})=\rho(T_{(k,k,k+1,k-1)})=\rho(T_{(k-1,k+1,k+1,k-1)}).$$
By Lemma \ref{lm2.3}, it is easy to see $$\rho(T_{(k,k,k,k)})=\rho(T_{(k-1,k)})=\rho(T_{(k-1,k+1,k+1,k-1)}).$$
Applying Lemma \ref{lm2.7} to these graphs, we get
\begin{eqnarray*}
\phi_{T_{(k,k,k,k)}}&=&p^{(k,k)}q^{(k)}x_2^{k-1}(x_2-x_1)\left(\frac{q^{(k,k)}}{p^{(k,k)}}-\frac{p^{(k)}}{q^{(k)}}x_1^{2k-2}\right),\\
\phi_{T_{(k,k,k+1,k-1)}}&=&p^{(k,k)}q^{(k-1)}x_2^k(x_2-x_1)\left(\frac{q^{(k,k)}}{p^{(k,k)}}-\frac{p^{(k-1)}}{q^{(k-1)}}x_1^{2k}\right),\\
\phi_{T_{(k,k,k+1,k-1)}}&=&p^{(k-1,k+1)}q^{(k)}x_2^{k-1}(x_2-x_1)\left(\frac{q^{(k-1,k+1)}}{p^{(k-1,k+1)}}-\frac{p^{(k)}}{q^{(k)}}x_1^{2k-2}\right),\\
\phi_{T_{(k-1,k+1,k+1,k-1)}}&=&p^{(k-1,k+1)}q^{(k-1)}x_2^k(x_2-x_1)\left(\frac{q^{(k-1,k+1)}}{p^{(k-1,k+1)}}-\frac{p^{(k-1)}}{q^{(k-1)}}x_1^{2k}\right).
\end{eqnarray*}

Let $\rho=\rho(T_{(k,k,k,k)})=\rho(T_{(k-1,k+1,k+1,k-1)})$ and $\rho'=\rho(T_{(k,k,k+1,k-1)})$. Write $J(\lambda)=p^{(k,k)}q^{(k-1)}x_2^k(x_2-x_1)$ and $K(\lambda)=p^{(k-1,k+1)}q^{(k)}x_2^{k-1}(x_2-x_1)$.
By Lemma \ref{lm3.2}, $J(\rho)>0$ and $K(\rho)>0$.

Note that $\rho$ is the root of both equations
\begin{equation}\label{eq6}
\frac{q^{(k,k)}}{p^{(k,k)}}=\frac{p^{(k)}}{q^{(k)}}x_1^{2k-2} \qquad \mbox{and} \qquad \frac{q^{(k-1,k+1)}}{p^{(k-1,k+1)}}=\frac{p^{(k-1)}}{q^{(k-1)}}x_1^{2k}.
\end{equation}

Note that  $\rho'$ is  the root of both equations
\begin{equation}\label{eq7}
\frac{q^{(k,k)}}{p^{(k,k)}}=\frac{p^{(k-1)}}{q^{(k-1)}}x_1^{2k} \qquad \mbox{and} \qquad \frac{q^{(k-1,k+1)}}{p^{(k-1,k+1)}}=\frac{p^{(k)}}{q^{(k)}}x_1^{2k-2}.
\end{equation}
We have 
\begin{eqnarray*}
\phi_{T_{(k,k,k+1,k-1)}}(\rho)&=&J(\rho)\left .\left(
\frac{p^{(k)}}{q^{(k)}}x_1^{2k-2}-\frac{p^{(k-1)}}{q^{(k-1)}}x_1^{2k}\right)\right |_{\rho}\\
&=&K(\rho)\left .\left(
\frac{p^{(k-1)}}{q^{(k-1)}}x_1^{2k}-\frac{p^{(k)}}{q^{(k)}}x_1^{2k-2}\right)\right |_{\rho}.
\end{eqnarray*}
Thus, $\phi_{T_{(k,k,k+1,k-1)}}(\rho)^2=-J(\rho)K(\rho)\left(x_1^{2k-2}\frac{p^{(k)}}{q^{(k)}}-\frac{p^{(k-1)}}{q^{(k-1)}}x_1^{2k}\right)^2\Big |_{\rho}\leq 0$. We get $\phi_{T_{(k,k,k+1,k-1)}}(\rho)=0$.
Similarly, we can prove $\phi_{T_{(k,k,k,k)}}(\rho')=0$. Hence, we get $\rho=\rho'$.

Now we prove $\rho(T_{(k,k,k+1,k-1)})<\rho(T_{(k,k+1,k,k-1)})$.
By Lemma \ref{lm2.8}, we have
$$\phi_{T_{(k,k,k+1,k-1)}}-\phi_{T_{(k,k+1,k,k-1)}}=(x_1-x_2)\left (p^{(k)}q^{(k-1)}-q^{(k)}p^{(k-1)}\right )=d_1d_2\lambda^2(\lambda^2-1)^2>0$$
for any $\lambda>\lambda_0$.
So $\rho(T_{(k,k,k+1,k-1)})<\rho(T_{(k,k+1,k,k-1)})$. We are done in this case.

\textbf{Case 2.} $\sum\limits_{i=1}^4 k_i=4k+1$.

Similarly, by Lemma \ref{lm2.8}, we have
\begin{eqnarray*}
\phi_{T_{(k,k+1,k,k)}}-\phi_{T_{(k,k+1,k+1,k-1)}}&=&(x_1-x_2)\left (p^{(k,k+1)}q_0-q^{(k,k+1)}p_0\right )\\
&=&\frac{(d_2x_1+1)(\lambda^2-1)^2x_2^{2k+1}}{(x_2-x_1)^3}\left (d_2^3-2d_1d_2x_1^{2k+1}-d_1^3x_1^{4k+2}\right ).
\end{eqnarray*}

Here we use proof by contradiction. Suppose $G_{n,n-8}^{min}=T_{(k,k+1,k+1,k-1)}$.
By Lemma \ref{lm3.5}, $d_2=\sqrt{\overline{c}d_1}x_1^{k}$ at $\lambda=\rho(T_{(k,k+1,k+1,k-1)})$.
Note $\overline{c}\to \lambda_0$ as $n \to \infty$.
When $n$ is large enough, we will get $\overline c > (2+\epsilon) x_1$ for some constant
$\epsilon>0$. Thus, we get
$$d_2^2= \overline{c}d_1 x_1^{2k}>(2+\epsilon)d_1x_1^{2k+1}.$$
For $n$ large enough, we have $\phi_{T_{(k,k+1,k,k)}}-\phi_{T_{(k,k+1,k+1,k-1)}}>0$
 at $\lambda=\rho(T_{(k,k+1,k+1,k-1)})$. Equivalently $\phi_{T_{(k,k+1,k,k)}}(\rho(T_{(k,k+1,k+1,k-1)}))>0$.
By Lemma \ref{lm3tree}, we get  $\rho(T_{k,k+1,k,k})<\rho(T_{k,k+1,k+1,k-1})$. Contradiction!
Hence, we have $G_{n,n-8}^{min}=T_{k,k+1,k,k}$.

\textbf{Case 3.} $\sum\limits_{i=1}^4 k_i=4k+2$. There is only one possible graph $T_{(k,k+1,k+1,k)}$.

\textbf{Case 4.} $\sum\limits_{i=1}^4 k_i=4k+3$.

Similarly by Lemma \ref{lm2.8}, we have
\begin{eqnarray*}
& &\phi_{T_{(k,k+1,k+1,k+1)}}-\phi_{T_{(k,k+1,k+2,k)}}\\
&=&(x_1-x_2)\left (p^{(k,k+1)}q_0-q^{(k,k+1)}p_0\right )\\
&=&\frac{(d_2x_1+1)(\lambda^2-1)^2x_2^{2k+1}}{(x_2-x_1)^3}\left (d_2^3-2d_1d_2x_1^{2k+1}-d_1^3x_1^{4k+2}\right )\\
&<&\frac{(d_2x_1+1)(\lambda^2-1)^2x_2^{2k+1}}{(x_2-x_1)^3}\left (d_2^3-2d_1d_2x_1^{2k+1}\right )\\
&=&\frac{d_2(d_2x_1+1)(\lambda^2-1)^2x_2^{2k+1}}{(x_2-x_1)^3}\left (d_2^2-2d_1x_1^{2k+1}\right ).
\end{eqnarray*}

We now suppose $G_{n,n-8}^{min}=T_{(k,k+1,k+1,k+1)}$ in this case. By Lemma \ref{lm3.5}, $d_2=\sqrt{\overline{c}d_1}x_1^{k+1}$ at $\lambda=\rho(T_{(k,k+1,k+1,k+1)})$.
Recall that $\overline{c}\to \lambda_0$ as $n \to \infty$.
When $n$ is large enough, we get $\overline c < 2x_2$.
Thus $d_2=\sqrt{\overline{c}d_1}x_1^{k+1}<\sqrt{2d_1x_2}x_1^{k+1}$.
We get $\phi_{T_{(k,k+1,k+2,k)}}(\rho(T_{(k,k+1,k+1,k+1)}))>0$.
Applying Lemma \ref{lm3tree} with $G_2=T_{(k,k+1,k+2,k)}$ and $G_1=T_{(k,k+1,k+1,k+1)}$,
we have $\rho(T_{k,k+1,k+2,k})<\rho(T_{k,k+1,k+1,k+1})$. Contradiction!
Hence $G_{n,n-8}^{min}=T_{k,k+1,k+2,k}$.
The proof is completed.
$\hfill\Box$

\end{document}